\documentclass[12pt]{article}

    \usepackage{psfrag}

    \usepackage{graphics}

    \usepackage{epsfig}
    \usepackage{amsthm,amsfonts,amsmath}
    \usepackage{color}   
    
    \usepackage[dvips]{geometry}
  \newcommand{\R}{\mathbb{R}}

\newtheorem{teorema}{Theorem}[section]
\newtheorem{lema}{Lemma}[section]
\newtheorem{prop}{Proposition} [section]

\newtheorem{obs}{Remark} [section]
\newtheorem{cor}{Corollary} [section]
\newtheorem{defi}{Definition} [section]

\newtheorem{example}{Example} [section]

\begin {document}

\begin{center}\Large{ \textbf{Rotational and Parabolic Surfaces in $\widetilde{PSL}_{2}(\R,\tau)$ and Applications}}
\end{center}

\vspace{0.2cm}

\begin{center}\large{By Carlos Espinoza Pe\~nafiel \footnote{The author was supported by FAPERJ - Brasil. \\ Departamento de Matematicas, Pontificia Universidade Catolica do Rio de Janeiro, RJ - Brasil.\\ e-mail: carloses@mat.puc-rio.br}\\ }
\end{center}

\vspace{0.8cm}

\textbf{Abstract.} We study surfaces of constant mean curvature which are invariant by one-parameter group of either rotational isometries or parabolic isometries, immersed into the homogeneous manifold $\widetilde{PSL}_2(\R,\tau)$. Also, we give some applications.
\vspace{0.5cm}

\vspace{0.5cm}

\textbf{Acknowledgements.} The author would like to thank his thesis advisors, Professor Ricardo Sa Earp and Professor Harold Rosenberg, for suggesting this problem, for many interesting and stimulating discussion on the subject and for the constant support throughout his work.
\vspace{0.5cm}

\section{Introduction}

\indent

In this paper we study constant mean curvature surfaces ($H$ surfaces), immersed into $\widetilde{PSL}_2(\R,\tau)$. In the papers Screw Motion Surfaces in $\mathbb{H}^2\times\R$ and $\mathbb{S}^2\times\R$, Ricardo Sa Earp and Eric Toubiana (see \cite{RT1}) studied the geometry of screw motions surfaces of constant mean curvature, in particular they studied the rotational surfaces, that is, when the pitch is zero. The study of the behavior of constant mean curvature surfaces invariant by rotational isometries immersed into $\mathbb{H}^2\times\R$ is given in the appendix of the paper Uniqueness of $H$-surfaces in $\mathbb{H}^2\times\R$, $|H|\leq1/2$, with boundary one or two parallel horizontal circles (see \cite{BRWE}).

In his Doctoral Thesis, Rami Younes (see \cite{RY}) gave several examples of minimal surfaces invariant by one-parameter group of isometries immersed into $\widetilde{PSL}_2(\R)$, there he explore the equation of the mean curvature in the divergence form to obtain the integral form for the function that determine the minimal surfaces.

We follow the ideas of \cite{RT1}, \cite{BRWE}, and \cite{RY} to obtain explicit formulas for the generating curve of the rotational and parabolic surfaces, see Lemma \ref{EFR} and Lemma \ref{EFP}. We give all rotational and parabolic $H$ surfaces in $\widetilde{PSL}_{2}(\R,\tau)$. For instance, we give a explicit $H=1/2$ rotational surface which is an entire graph and it is asymptotic to the asymptotic boundary of $\widetilde{PSL}_2(\R,\tau)$. This surface has the following expression, (see Example \ref{ER1})
\begin{equation*}
   u(\rho)=2\sqrt{\cosh(\rho)}-2\arctan(\sqrt{\cosh(\rho)})
\end{equation*}

We generalize part of the works \cite{RT1}, \cite{BRWE}, since in the case $\tau=0$, we recover the $H$ surfaces invariant by rotational isometries in $\mathbb{H}^2\times\R$.

Also, we focus our attention on $H$ surfaces invariant by parabolic isometries immersed into $\widetilde{PSL}_2(\R,\tau)$. Here we give all such surfaces. Again, we generalize part of the work Parabolic and Hyperbolic Screw Motion Surfaces in $\mathbb{H}^{2}\times\R$ (see \cite{R}), that is, if $\tau\equiv0$ we recover the $H$ surfaces invariant by parabolic isometries in $\mathbb{H}^2\times\R$.

The geometry of rotational and parabolic isometries in $\widetilde{PSL}_2(\R,\tau)$ is analogous to the geometry of rotational and parabolic isometries in $\mathbb{H}^2\times\R$, see \cite{RT1}, \cite{BRWE}, and \cite{R} respectively. For example, there is also a notion of growth for rotational surface in $\widetilde{PSL}_2(\R,\tau)$ as well as $\mathbb{H}^2\times\R$, see \cite{BS}. The growth for rotational surfaces in $\widetilde{PSL}_2(\R,\tau)$  is proportional to the growth of rotational surfaces in $\mathbb{H}^2\times\R$ (see Applications).

We give some applications, for instance, we prove that, there is no entire graph in $\widetilde{PSL}_2(\R,\tau)$ with constant mean curvature $H>1/2$ (see Proposition \ref{AP1}), which is a step to prove a more general theorem, see Theorem \ref{MT1}. Also, we give a halspace theorem for $H=1/2$ surface in $\widetilde{PSL}_2(\R,\tau)$, see Theorem \ref{THS1}, this theorem was proved by Barbara Nelly and Ricardo Sa Earp in \cite[Theorem 1]{BS}, when the ambient space is $\mathbb{H}^2\times\R$. The proof in the space $\widetilde{PSL}_2(\R,\tau)$ is analogous.

The details of the proof of the above result as well as the study of surfaces of constant mean curvature invariant by hyperbolic isometries immersed into $\widetilde{PSL}_2(\R,\tau)$ and other applications are given in \cite{CP}.

\section{Notations}

There are different models for the Hyperbolic space; in this paper we work with the half-plane model and the disk model, in both cases the Hyperbolic space will be denote by $M^2$. More precisely:\\
Let $M^2$ be the self-plane model of the Hyperbolic space, so
\begin{equation*}
    M^{2}:=\mathbb{H}^{2}=\{(x,y)\in\R^{2},y>0\}
\end{equation*}
endowed with metric
\begin{equation*}
    d\sigma^{2}=\lambda^{2}(dx^{2}+dy^{2}), \hspace{0.5cm} \lambda=\frac{1}{y}
\end{equation*}
Let $M^2$ be the disk model of the Hyperbolic space, so
\begin{equation*}
    M^{2}:=\mathbb{D}^{2}=\{(x,y)\in\R^{2},x^{2}+y^{2}<1\}
\end{equation*}
endowed with metric
\begin{equation*}
    d\sigma^{2}=\lambda^{2}(dx^{2}+dy^{2}), \hspace{0.5cm} \lambda=\frac{2}{1-(x^{2}+y^{2})}
\end{equation*}

The natural orthonormal frame field on $M^{2}$ is given by $\{e_{1},e_{2}\}$, where,

\begin{equation*}
 \left\{
        \begin{array}{ll}
          e_{1} = \lambda^{-1}\partial_{x}; & \hbox{  $$} \vspace{0.5cm} \\
          e_{2} = \lambda^{-1}\partial_{y}; & \hbox{ $$}
        \end{array}
      \right.
\end{equation*}
Let $\{E_1,E_2,E_3\}$ be the orthormal frame of $\widetilde{PSL}_{2}(\R)$ where $\{E_1,E_2\}$ is the horizontal lift of the frame $\{e_1,e_2\}$, and $E_3$ is the Killing field tangent to the fibers. We say that $\{E_1,E_2\}$ are horizontal field and $E_3$ is a vertical field, see Lemma \ref{F1}.

\section{ The Space $\widetilde{PSL}_{2}(\R,\tau)$}

It is natural to study $H$ surfaces invariant by isometries immersed into $\widetilde{PSL}_2(\R,\tau)$, since both $\mathbb{H}^2\times\R$ and $\widetilde{PSL}_2(\R,\tau)$ are two of the eight Thurston's geometries and there exists a natural projection (see \cite{WT})
\begin{equation*}
    \pi:E\longrightarrow M^2, \hspace{0.5cm} \pi(x,y,t)=(x,y)
\end{equation*}
where either $E=\widetilde{PSL}_2(\R,\tau)$ or $E=\mathbb{H}\times\R$ and $M^2$ is the Hyperbolic space, this projection $\pi$ is a Riemannian submersion.

The space $\widetilde{PSL}_{2}(\R)$ is a simply connected homogeneous manifold whose isometry groups has dimension 4, such a manifold is a Riemannian fibration over the 2-dimensional Hyperbolic space, the fibers are geodesics and there exists a one-parameter family of translations along the fibers, generated by a unit Killing field $E_{3}$, which will be called the vertical vector field.

To construct the space $\widetilde{PSL}_{2}(\R,\tau)$ we follow ideas of the professor Eric Toubiana (see \cite{T}). The simply connected homogeneous 3-manifolds $\widetilde{PSL}_{2}(\R,\tau)$ is given by,
\begin{equation*}
    \widetilde{PSL}_{2}(\R,\tau)=\{(x,y,t)\in\R^3;(x,y)\in M^2,t\in\R\}
\end{equation*}
endowed with metric,
\begin{equation*}
    g:=ds^2=\lambda^2(dx^2+dy^2)+(2\tau(\dfrac{\lambda_y}{\lambda}dx-\dfrac{\lambda_x}{\lambda} dy)+dt)^2, \hspace{0.5cm} \lambda=\dfrac{2}{1-(x^2+y^2)}
\end{equation*}
By considering the Riemannian submersion $\pi:\widetilde{PSL}_{2}(\R,\tau)\longrightarrow M^2$, we obtain the next Lemma.

\begin{lema}\label{F1} The fields $\{E_{1},E_{2},E_{3}\}$ in the referential $\{\partial_{x},\partial_{y},\partial_{t}\}$ are given by,

\begin{equation*}
     E_{1}=\frac{1}{\lambda}\partial_{x}-2\tau\frac{\lambda_{y}}{\lambda^{2}}\partial_{t}, \hspace{1cm} E_{2}=\frac{1}{\lambda}\partial_{y}+2\tau\frac{\lambda_{x}}{\lambda^{2}}\partial_{t}, \hspace{1cm} E_{3}=\partial_{t}.
\end{equation*}
\end{lema}
\begin{proof} Use that $\{E_1,E_2\}$ is the horizontal lift of $\{e_1,e_2\}$ and $\{E_1,E_2,E_3\}$ is an orthonormal frame. (See \cite{T})
\end{proof}

\subsection{Isometries in $\widetilde{PSL}_{2}(\R,\tau)$}

Since there exist a Riemannian submersion $\pi:\widetilde{PSL}_{2}(\R,\tau)\longrightarrow M^2$, the isometries of $\widetilde{PSL}_{2}(\R,\tau)$ are strongly related with the isometries of the Hyperbolic space $M^2$. It is well know the behavior of the isometries on $M^2$ (see \cite{RT}), so we can know the behavior of the isometries on $\widetilde{PSL}_{2}(\R,\tau)$, see Proposition \ref{i01}.\newline

From now on we identify the Euclidean space $\R^2$ with the set of complex numbers $\mathbb{C}$, more precisely $z=x+iy\approx(x,y)$ . So, if we take $M^2\equiv\mathbb{D}^2$, we obtain
\begin{equation*}
    \widetilde{PSL}_{2}(\R)=\{(z,t)\in\R^3;x^2+y^2<1,t\in\R\}
\end{equation*}
endowed with metric,

\begin{equation*}
    d\sigma^2=\lambda^{2}(z)|dz|^{2}+(i\tau\lambda(\overline{z}dz-zd\overline{z})+dt)^{2}
\end{equation*}
If we take $M^2\equiv\mathbb{H}^2$, we obtain
\begin{equation*}
    \widetilde{PSL}_{2}(\R)=\{(z,t)\in\R^3;y>0,t\in\R\}
\end{equation*}
endowed with metric,

\begin{equation*}
    d\sigma^2=\lambda^{2}(z)|dz|^{2}+(-\tau\lambda(dz+d\overline{z})+dt)^{2}
\end{equation*}

Let $F$ an isometry of $\widetilde{PSL}_{2}(\R,\tau)$, since $\pi:\widetilde{PSL}_{2}(\R,\tau)\longrightarrow M^{2}$ is a Riemannian submersion, then the horizontal component $\pi\circ F$ is an isometry of $M^{2}$, so $F(z,t)=(f(z),h(z,t))$. where $f$ is an isometry of $M^{2}$.

\begin{prop} \label{i01} The isometries positives of $\widetilde{PSL}_{2}(\R,\tau)$ are given by,
\begin{equation*}
    F(z,t)=(f(z),t-2\tau\arg f^{\prime}+c)
\end{equation*}
where $f$ is a positive isometry of $M^{2}$ and $c$ is a real number.
\end{prop}
\begin{proof} This follow after a long computation by considering that $F$ preserve the metric (see \cite{T}).
\end{proof}

\subsection{Graph in $\widetilde{PSL}_{2}(\R,\tau)$}

Since, $\pi:\widetilde{PSL}_{2}(\R,\tau)\longrightarrow M^2$ is a Riemannian submersion, this is possible to speak of graph in $\widetilde{PSL}_{2}(\R,\tau)$.

\begin{defi} A graph in $\widetilde{PSL}_{2}(\R,\tau)$ is the image of a section $s_0:\Omega\subset M^2\longrightarrow\widetilde{PSL}_{2}(\R,\tau)$, where $\Omega$ is a domain of $M^2$.
\end{defi}

Given a domain $\Omega\subset M^{2}$ we also denote by $\Omega$ its lift to $M^{2}\times\{0\}$, with this identification we have that the graph $\Sigma (u)$ of $u\in (C^{0}(\partial\Omega)\cap C^{\infty}(\Omega))$ is given by,
 \begin{equation*}
    \Sigma(u)=\{(x,y,u(x,y))\in\widetilde{PSL}_{2}(\R,\tau); (x,y)\in\Omega\}
 \end{equation*}
\begin{lema}\label{l1} Let $\Sigma(u)$ be the graph in $\widetilde{PSL}_{2}(\R,\tau)$ of the function $u:\Omega\subset M^{2}\longrightarrow\R$ with constant mean curvature $H$. Then, the function $u$ satisfies the equation
 \begin{equation*}
    2H=div_{M^{2}}\displaystyle(\frac{\alpha}{W}e_{1}+\frac{\beta}{W} e_{2}),
 \end{equation*}
 where $W=\displaystyle\sqrt{1+\alpha^{2}+\beta^{2}}$ and,

\begin{itemize}
  \item $\alpha=\displaystyle\frac{u_{x}}{\lambda}+2\tau\frac{\lambda_{y}}{\lambda^{2}}$,
  \item $\beta=\displaystyle\frac{u_{y}}{\lambda}-2\tau\frac{\lambda_{x}}{\lambda^{2}}$.
\end{itemize}

\end{lema}
\begin{proof} The proof follow since $E_3$ is a Killing field and $\pi:\widetilde{PSL}_{2}(\R,\tau)\longrightarrow M^2$ is a Riemannian submersion. See \cite{RY}.
\end{proof}

\subsection{The Mean Curvature Equation in $\widetilde{PSL}_{2}(\R,\tau)$}

By considering the coefficients of the first and second fundamental form, of a surface immersed into $\widetilde{PSL}_{2}(\R,\tau)$ we can write the mean curvature equation, for example taking the half plane model $M^2$ for the hyperbolic space, the mean curvature equation for a vertical graph is given by:
\begin{equation*}
    2H\lambda^2 m^3=u_{xx}(\lambda^3+\lambda u^{2}_y)+u_{yy}\lambda(\lambda^2+(u_x-2\tau\lambda)^2) -2u_{xy}\lambda(u_x-2\tau\lambda)u_y-u_xu_y\lambda^2(u_x-2\tau\lambda)-\lambda^2u^3_y
\end{equation*}
where $m=\sqrt{\lambda^2+(2\tau\lambda-u_x)^2+u_y}$ (see \cite{CP}).

\subsection{Maximum Principle in $\widetilde{PSL}_{2}(\R,\tau)$}

An important criterium in Riemannian Geometry is the maximum principle. \\ \\

$\mathbf{Maximum}$ $\mathbf{Principle.}$ Let $S_1$ and $S_2$ be two surfaces with constant mean curvature $H$ that are tangent at a point $p\in int(S_1)\cap int(S_2)$. Assume that the mean curvature vectors of $S_1$ and $S_2$ at $p$ coincide and that, around $p$, $S_1$ lies on one side of $S_2$. Then $s_2\equiv S_2$. When the intersection point $p$ belongs to the boundary of the surfaces, the result holds as well, provided further that the two boundaries are tangent and both are local graphs over a common neighborhood in $T_{p}S_1=T_{p}S_2$.

\section{Screw Motions Surfaces in $\widetilde{PSL}_{2}(\R,\tau)$}

A screw motion surface, is a surface which is invariant by one-parameter group of isometries in $\widetilde{PSL}_{2}(\R,\tau)$, this group of isometries is the composition of rotational isometries (prevenient of rotational isometries of $\mathbb{D}^2$ ), together with vertical translations in a proportional way.

We focus our attention on rotational surfaces, but for later use, we give the integral form for a screw motion surface. The idea is simple, we will take a curve in the $xt$ plane and we will apply one-parameter group of Rotational isometries together vertical translation to obtain a screw motion surface on $\widetilde{PSL}_{2}(\R,\tau)$. We denote by $\alpha(x)=(x,0,u(x))$ the curve in the $xt$ plane and by $S$ the screw motion surface generate by $\alpha$.

Since the most simple rotational isometry in $\mathbb{D}^{2}$ is the rotation around the origin, we re-parameterized the Hyperbolic disk whit coordinates $\rho$ and $\theta$, so that
\begin{eqnarray*}
  x &=& \tanh(\frac{\rho}{2})\cos(\theta) \\
  y &=& \tanh(\frac{\rho}{2})\sin(\theta)
\end{eqnarray*}
where $\rho$ is the hyperbolic distance measure from the origin of $\mathbb{D}^{2}$.

So, the surface $S$ is parameterized by,
 \begin{equation*}
    \varphi(\rho,\theta)=(\tanh(\frac{\rho}{2})\cos(\theta),\tanh(\frac{\rho}{2})\sin(\theta),u(\rho)- 2\tau\theta+\widetilde{l}\theta)
\end{equation*}
The next Lemma is crucial for our study. We follow ideas of Appendix A of \cite{LMH}. Denoting by $l=\widetilde{l}-2\tau$, we have the next lemma.
\begin{lema}\label{r1} With the notations above, and denoting by $H$ the mean curvature of $S$, then the function $u$ satisfies
\begin{equation*}
    u(\rho)=\displaystyle\int\frac{(2H\cosh(\rho)+d)\sqrt{1+ \left[\dfrac{l}{\sinh(\rho)}-2\tau\tanh\left(\dfrac{\rho}{2}\right)\right]^2}}{\sqrt{\sinh^{2}(\rho)-(2H\cosh(\rho)+c)^{2}}}
\end{equation*}
where $d$ is a real number, and $l$ is the pitch.
\end{lema}
\begin{proof}
Since $S$ has mean curvature $H$, then by lemma \ref{l1} the function $u$ satisfies the equation
\begin{equation}\label{d1}
    2H=div_{\mathbb{D}^{2}}\displaystyle(\frac{\alpha}{W}e_{1}+\frac{\beta}{W} e_{2}),
 \end{equation}
where $W=\displaystyle\sqrt{1+\alpha^{2}+\beta^{2}}$ and,

\begin{itemize}
  \item $\alpha=\dfrac{u_{x}}{\lambda}+2\tau y$,
  \item $\beta=\dfrac{u_{y}}{\lambda}-2\tau x$.
\end{itemize}

Set $X_{u}=\dfrac{\alpha}{W}e_{1}+\dfrac{\beta}{W} e_{2}$. then in coordinates $\rho$ and $\theta$, $X_{u}$ is given by
\begin{equation*}
    X_{u}=\dfrac{1}{W}\left[u_{\rho}\partial_{\rho}+ \left[\dfrac{l}{\sinh^2(\rho)}-2\tau\frac{\tanh(\rho/2)}{\sinh(\rho)}\right]\partial_{\theta} \right]
\end{equation*}
and
\begin{equation*}
    W=\displaystyle\sqrt{1+\left[\dfrac{l}{\sinh^2(\rho)}-2\tau\frac{\tanh(\rho/2)}{\sinh(\rho)}\right]^2+u^{2}_{\rho}}
\end{equation*}
Let $\theta_{0},\theta_{1}\in(0,2\pi)$ with $\theta_{0}<\theta_{1}$ and $\rho_{0},\rho_{1}\in\R$ with $\rho_{0}<\rho_{1}$ and consider the domain $\Omega=[\theta_{0},\theta_{1}]\times[\rho_{0},\rho_{1}]$. By integrating the equation \ref{d1}, we obtain

\begin{equation*}
    \displaystyle\int_{\partial(\Omega)} <X_{u},\eta> =2HArea([\theta_{0},\theta_{1}]\times[\rho_{0},\rho_{1}])
\end{equation*}
where $\eta$ is the outer co-normal. This gives,
\begin{equation*}
    \partial_{\rho}\left(\frac{u_{\rho}\sinh(\rho)}{W}\right)=2H\sinh(\rho)
\end{equation*}
by integrating this expression we get the lemma.

\end{proof}

\section{Rotational Surface in $\widetilde{PSL}_{2}(\R,\tau)$}

From Proposition \ref{i01}, we know that to obtain a rotational motion on $\widetilde{PSL}_{2}(\R,\tau)$ is necessary consider a rotational motion on $\mathbb{D}^2$, but by considering only a rotational motion on $\mathbb{D}^2$, the induced isometry on $\widetilde{PSL}_{2}(\R,\tau)$ gives a screw motion, since the vertical translations  are isometries on $\widetilde{PSL}_{2}(\R,\tau)$, we consider our one-parameter group $\Gamma$ of isometries as being the composition of rotational motion from $\mathbb{D}^2$ together vertical translations in such way that $\Gamma$ give exactly a Rotational motion on $\widetilde{PSL}_{2}(\R,\tau)$.

Our idea is simple, we take a curve in the $xt$ plane and we will apply one-parameter group $\Gamma$ of Rotational isometries on $\widetilde{PSL}_{2}(\R,\tau)$ to the curve to generate a rotational surface.

An immediate consequence of the Lemma \ref{r1} is the next corollary.
\begin{cor}\label{CR1}
Consider the graph $t=u(\rho)$ in the plane $xt$, and denote by $S=grap(\Gamma u)$, so the function $u$ satisfies the next equation
\begin{equation}\label{FR}
    u(\rho)=\displaystyle\int\frac{(2H\cosh(\rho)+d)\sqrt{1+ 4\tau^2\tanh^2\left(\dfrac{\rho}{2}\right)}}{\sqrt{\sinh^{2}(\rho)-(2H\cosh(\rho)+c)^{2}}}
\end{equation}
where $d$ is a real number.
\end{cor}

\begin{obs} There is other form to compute the equation of the mean curvature, by considering the first and second fundamental form (see section 3.3), a hard computation give the same result.
\end{obs}

Now, we explore the Corollary \ref{CR1}, by considering $\tau=-1/2$ we obtain the next consequences:
\begin{lema}\label{EFR} Setting $d=-2H$, then the integral
\begin{equation*}
    u(\rho)=\displaystyle\int\frac{(2H\cosh(\rho)-2H)\sqrt{1+ \tanh^{2}(\frac{\rho}{2})}}{\sqrt{\sinh^{2}(\rho)-(2H\cosh(\rho)-2H)^{2}}}d\rho
\end{equation*}
has the following solution,
\begin{itemize}
  \item If $4H^{2}-1>0$ then,
  \begin{equation*}
    u(\rho)=\frac{4\sqrt{2}H}{\sqrt{4H^{2}-1}}\left[\arctan\left(\frac{\sqrt{\cosh(\rho)}}{\frac{4H^{2}+1}{4H^{2}-1} -\cosh(\rho)}\right) \right] -2\arctan\left(\frac{\sqrt{\frac{8H^{2}}{4H^{2}-1}}\sqrt{\cosh(\rho)}}{\sqrt{\frac{4H^{2}+1}{4H^{2}-1} -\cosh(\rho)}}\right)
  \end{equation*}

  \item If $1-4H^{2}>0$ then,
 \begin{eqnarray*}
   u(\rho) &=& \displaystyle\frac{4\sqrt{2}H}{\sqrt{1-4H^{2}}} \ln\left(\sqrt{\cosh(\rho)}+\sqrt{\frac{1+4H^{2}}{1-4H^{2}}+\cosh(\rho)}\right)+ \\
     &+& 2\arctan\left(-\sqrt{\frac{8H^{2}}{1-4H^{2}}}\frac{\sqrt{\cosh(\rho)}} {\sqrt{\frac{1+4H^{2}}{1-4H^{2}}+\cosh(\rho)}}\right)
 \end{eqnarray*}

\end{itemize}
\end{lema}
The proof of the Lemma \ref{EFR} is a straightforward computation.

\begin{example} Putting $H=\sqrt{3}/2$, we obtain a rotational surface (noncomplete), which is a graph over a domain in $\mathbb{D}^2$, since the rotation by $\pi$ around the $x$ axis is a isometry of $\widetilde{PSL}_2(\R,\tau)$, the surface is actually a complete embedded rotational surface. By using Maple the graph is given by,

\begin{equation*}
    u(\rho)=2\sqrt{3}\arcsin\left(\dfrac{\sqrt{\cosh(\rho)}}{\sqrt{2}}\right) +2\arctan\left(\dfrac{-\sqrt{3}\sqrt{\cosh(\rho)}}{\sqrt{2-\cosh(\rho)}}\right)
\end{equation*}

\begin{center}
  \begin{tabular}{cc}
    \includegraphics[width=2.5in]{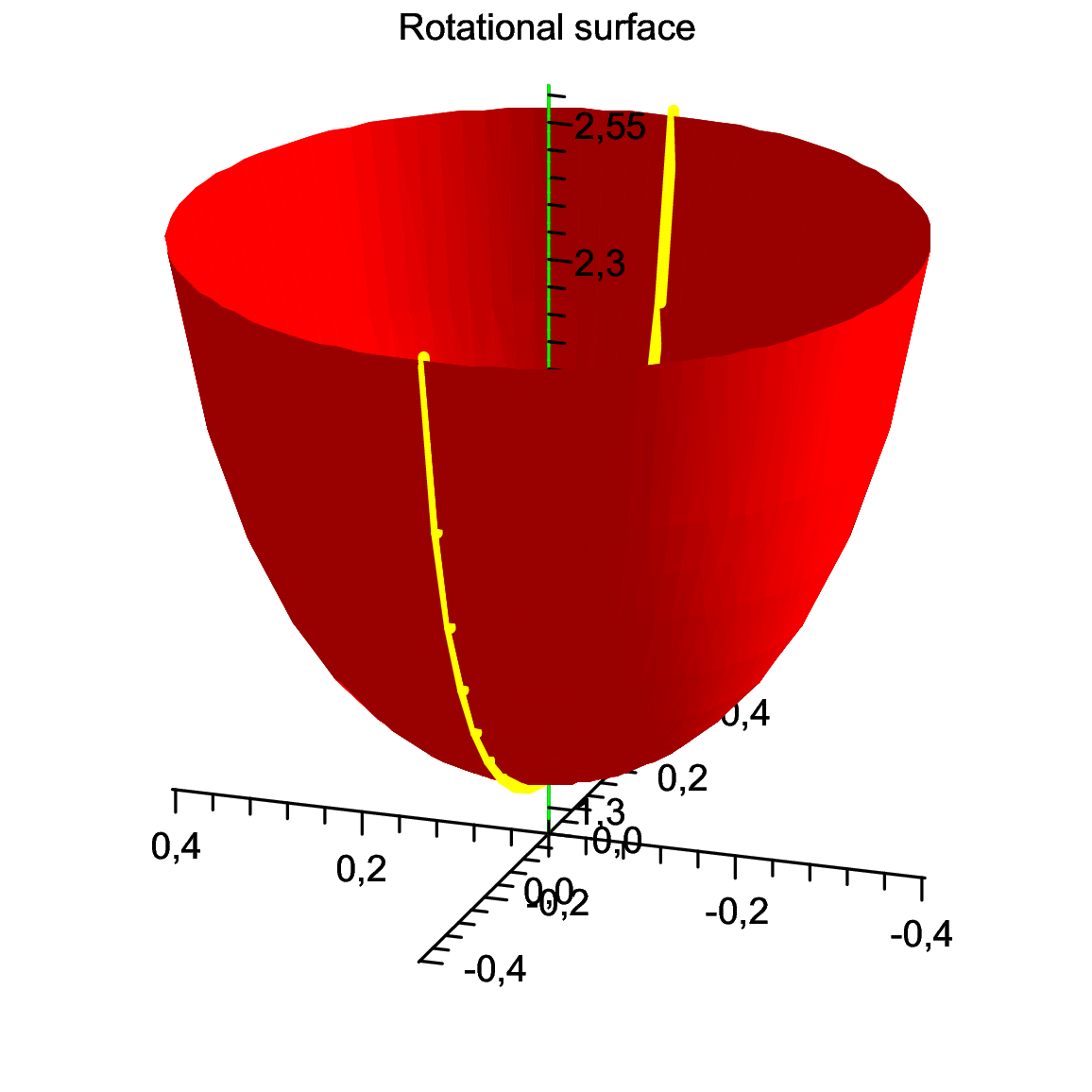} &
    \includegraphics[width=2.5in]{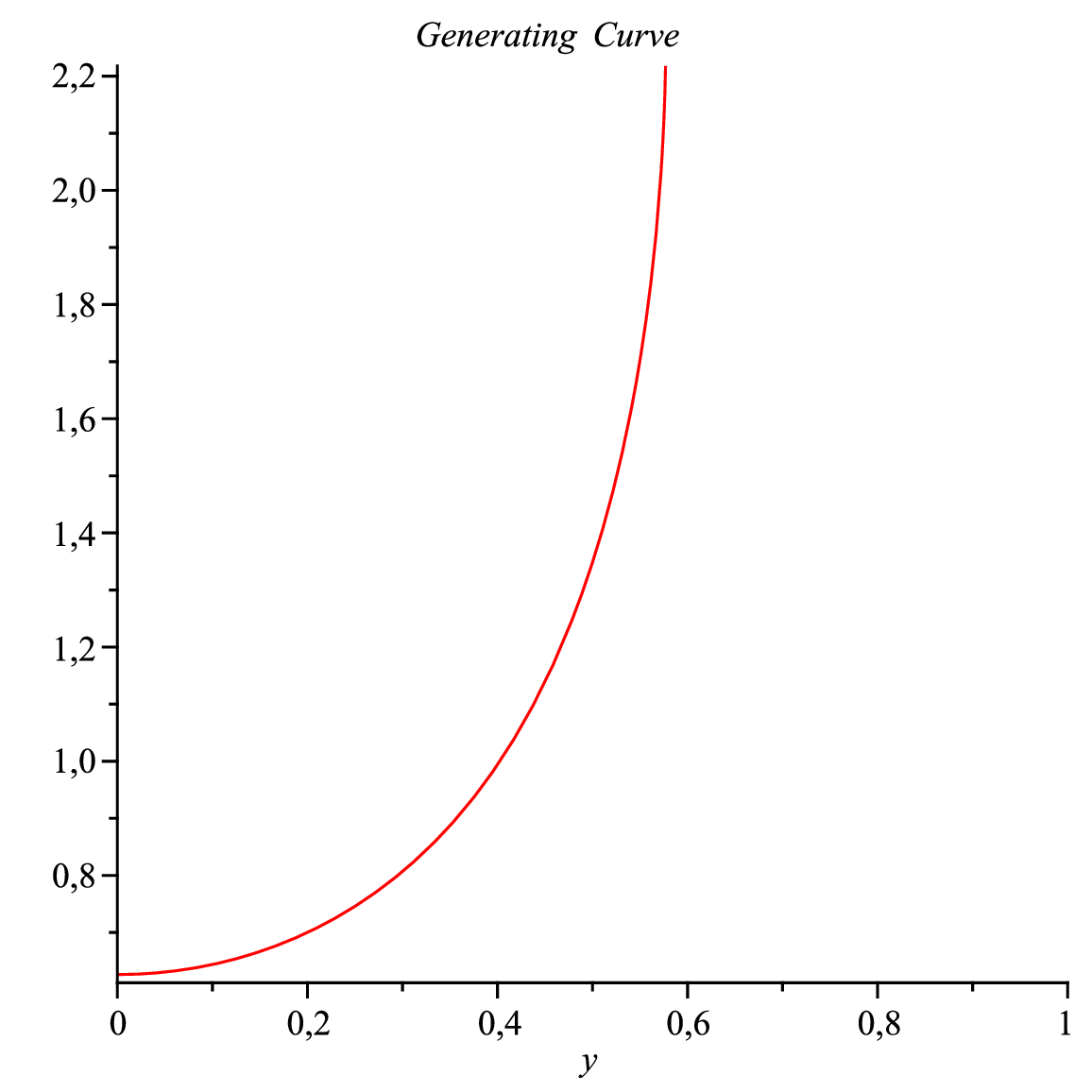}
  \end{tabular}
\end{center}

\end{example}

\begin{example}\label{ER1} Putting $d=-2H$ and $H=1/2$, then by integrating the formula \ref{FR}, we obtain $H=1/2$ surfaces invariant by rotations in $\widetilde{PSL}_{2}(\R)$ which is an entire graph, this surface is asymptotic to asymptotic boundary. More specifically
\begin{equation*}
   u(\rho)=2\sqrt{\cosh(\rho)}-2\tan^{-1}(\sqrt{\cosh(\rho)})
\end{equation*}
which expressed in Euclidian coordinates gives
\begin{equation*}
    u(x,y)=2\sqrt{\cosh(2\tanh^{-1}(\sqrt{x^{2}+y^{2}}))}-2\tan^{-1}(\sqrt{\cosh(2\tanh^{-1}(\sqrt{x^{2}+y^{2}}))})
\end{equation*}
Maple gives,

\begin{center}
  \begin{tabular}{cc}
    \includegraphics[width=2.5in]{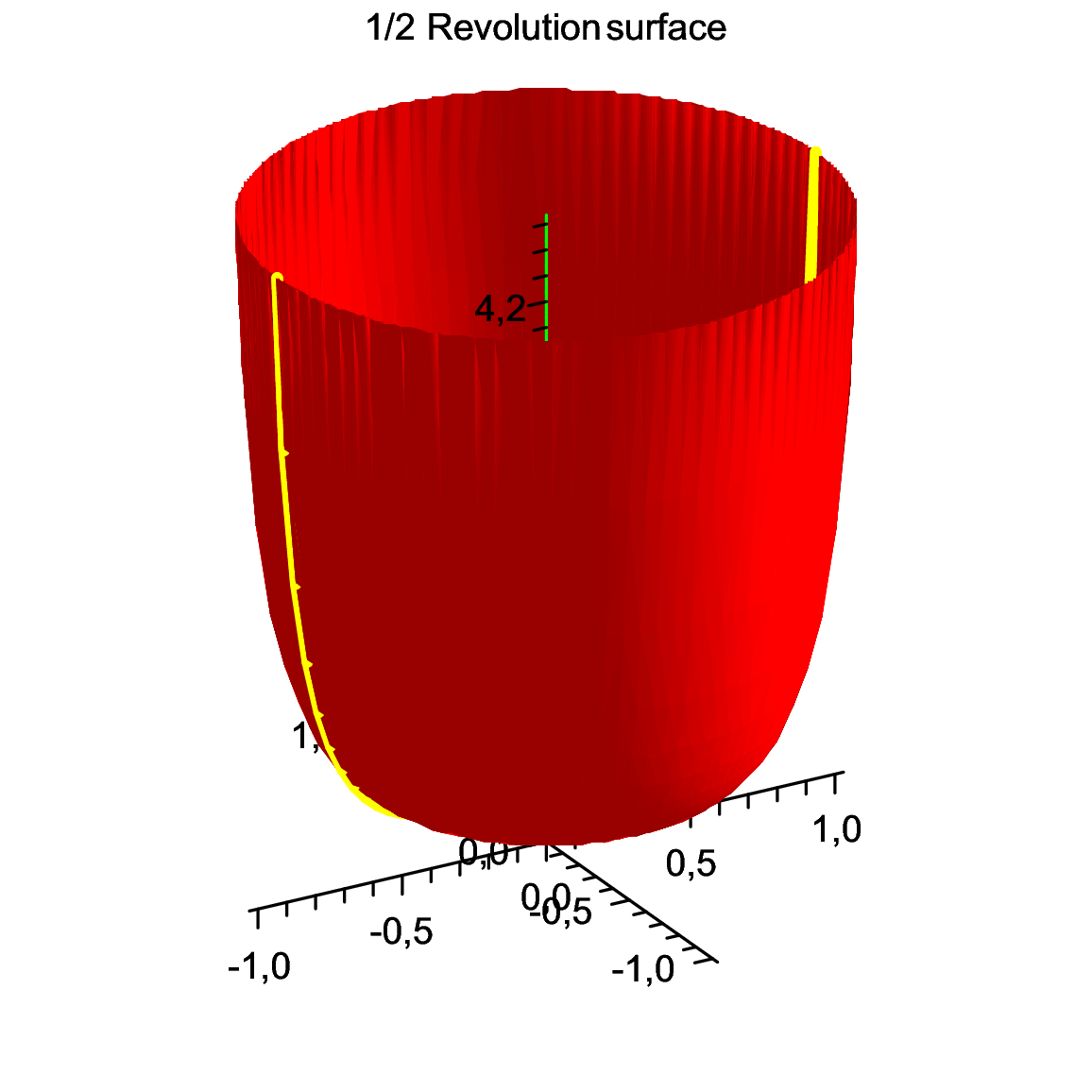} &
    \includegraphics[width=2.5in]{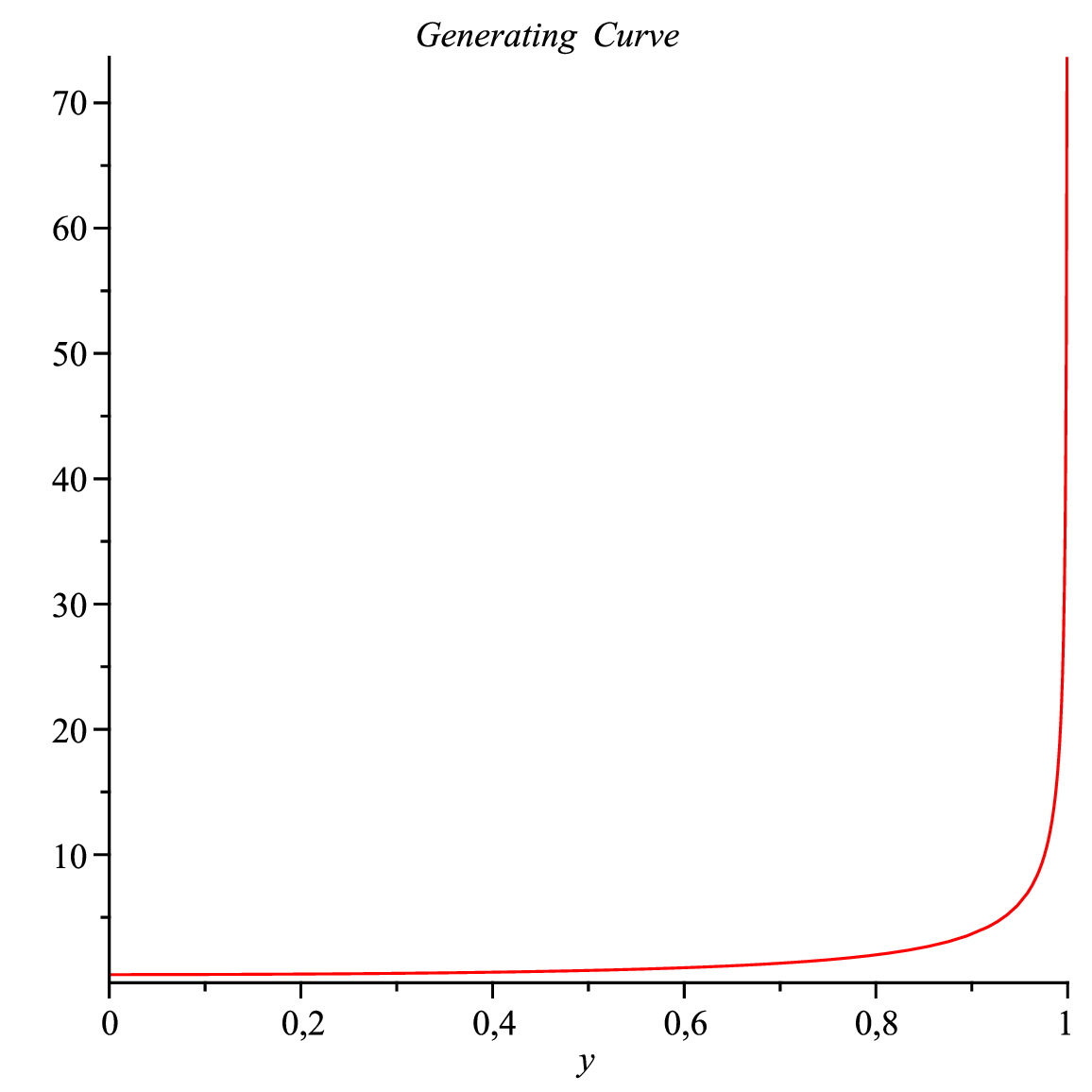}
  \end{tabular}
\end{center}

\end{example}

\section{Minimal Surfaces Invariant by Rotations in $\widetilde{PSL}_2(\R,\tau)$}

In this section we study quickly the behavior of minimal rotational surface, that is $H\equiv0$.

Rami Younes give a first integral for minimal rotational surfaces in $\widetilde{PSL}_2(\R)$ (see \cite{RY}). He gives examples of rotational minimal surfaces as well as hyperbolic and parabolic minimal surfaces. \\
By considering $H\equiv0$ in the Corollary \ref{CR1} we obtain the next proposition.

\begin{prop} (Minimal Rotational Surfaces) For each $d\geq0$ there exist a complete minimal rotational surface $\mathcal{M}_d$. The surface $\mathcal{M}_0$ is the slice $t=0$. For $d>0$ the rotational surface $\mathcal{M}_d$ (called catenoid) is embedded and homeomorphic  to an annulus.
\end{prop}
\begin{proof} Observe that, the Corollary \ref{CR1} gives,
\begin{equation*}
    u(\rho)=\displaystyle\int_{arcsinh(d)}^{\rho}\dfrac{d\sqrt{1+ 4\tau^2\tanh^{2}(\dfrac{r}{2})}}{\sqrt{\sinh^{2}(r)-d^{2}}}dr
\end{equation*}
A simple computation gives $u^{\prime}=\dfrac{d\sqrt{1+ 4\tau^2\tanh^{2}(\dfrac{r}{2})}}{\sqrt{\sinh^{2}(r)-d^{2}}}>0$ and $u^{\prime\prime}<0$

\end{proof}
\begin{example} With Maple's help, we plot the catenoid $\mathcal{M}_1$. Observe that by considering the rotation by $\pi$ around the $x$ axis, we obtain a complete embedded surface.
\begin{center}
  \begin{tabular}{cc}
    \includegraphics[width=2.5in]{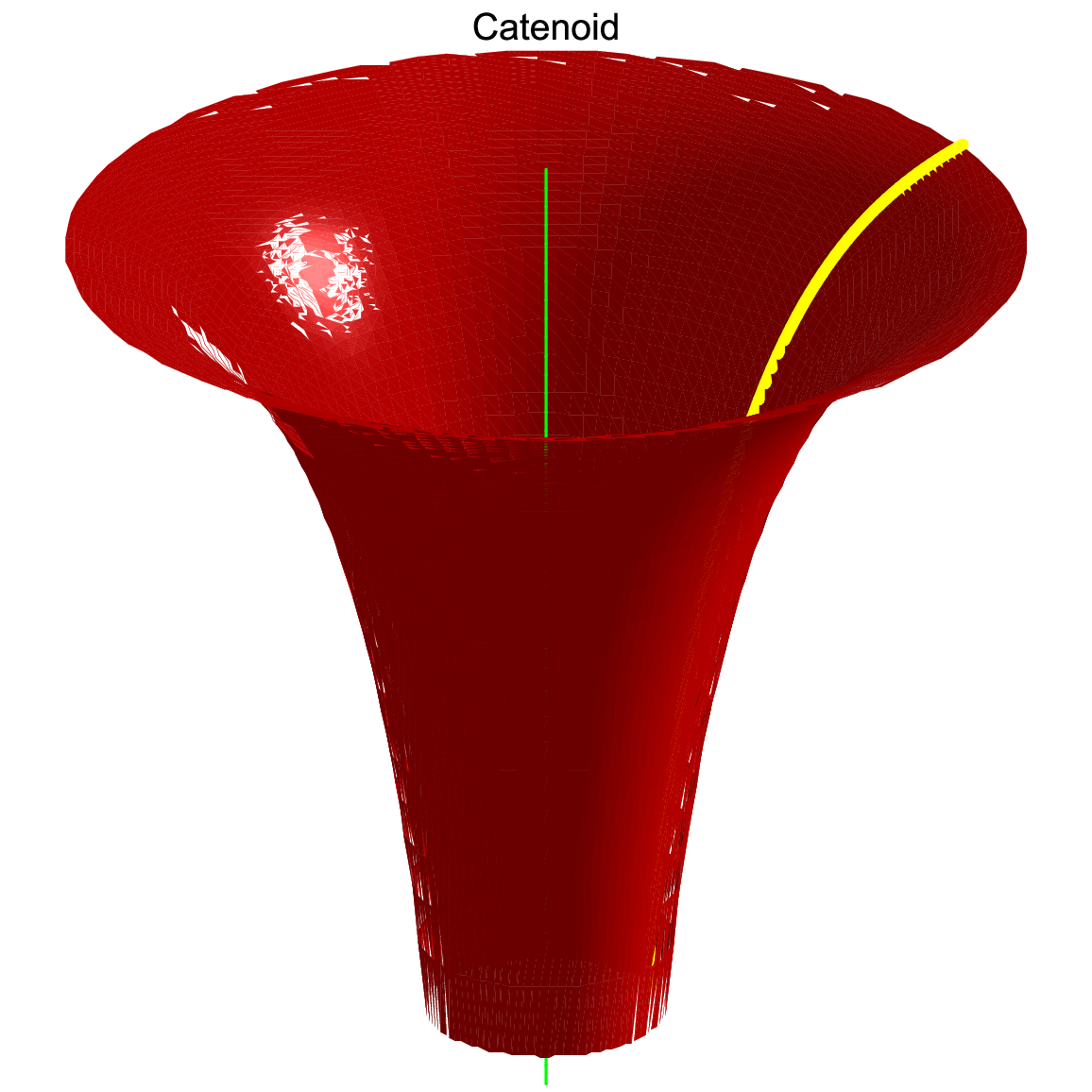} &
    \includegraphics[width=2.5in]{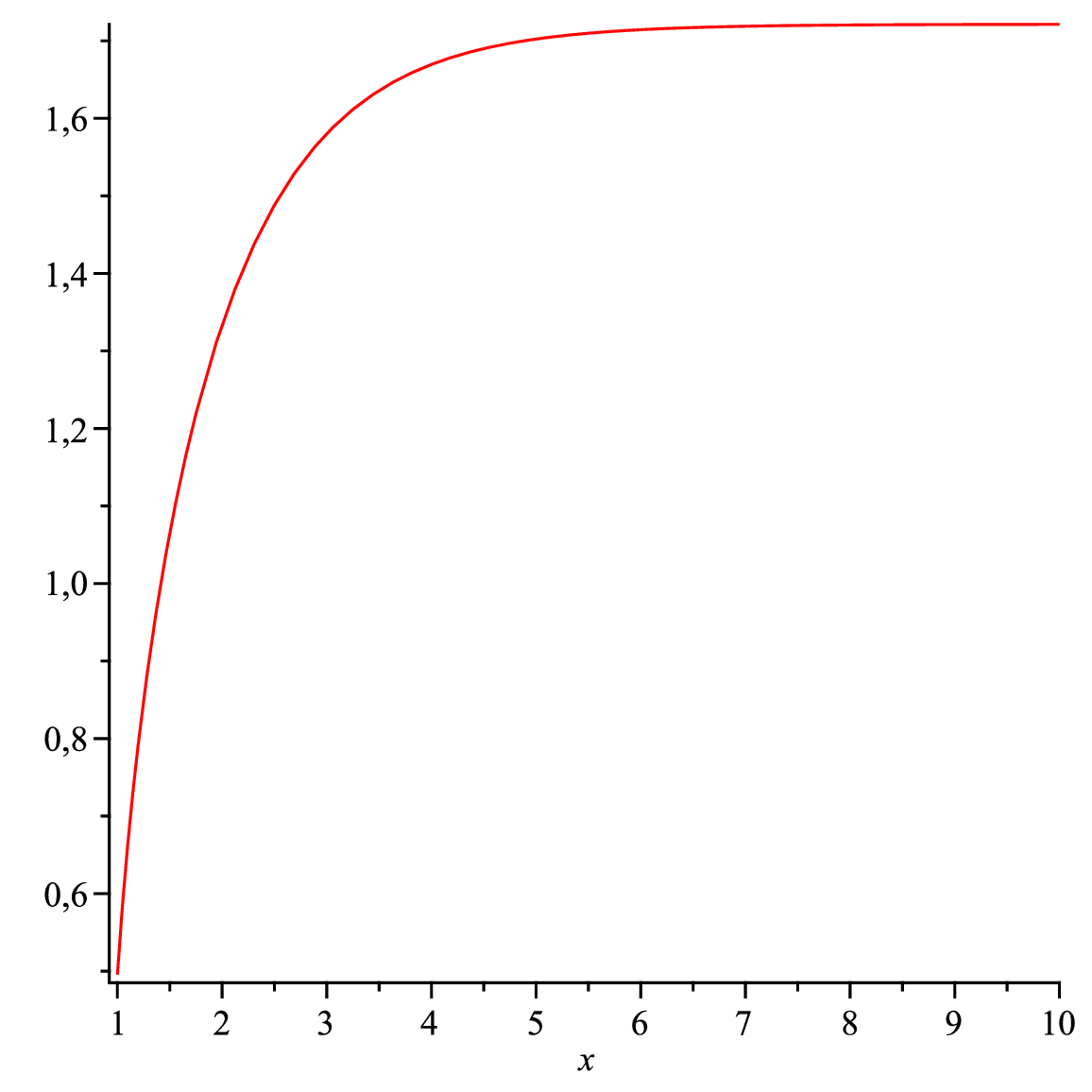}
  \end{tabular}
\end{center}
\end{example}

\section{Surfaces Invariant by Rotations in $\widetilde{PSL}_2(\R,\tau)$ with Constant Mean Curvature $H\neq0$}

In this section, we follow the ideas of the paper \cite[Proposition 5.2,Proposition 5.3]{BRWE}, to describe the behavior of rotational $H$-surfaces. For later use we define the functions $g(\rho)$ and $f(\rho)$ setting for $d\in\R$ and $H>0$.
\begin{eqnarray*}
  g(\rho) &=& d+2H\cosh(\rho) \\
  f(\rho) &=& \sinh^2(\rho)-(d+2H\cosh(\rho))^2 \\
    &=& (1-4H^2)\cosh^2(\rho)-1-d^2
\end{eqnarray*}

\begin{lema}\label{rl1} Assume $0<H<1/2$. We have $f(\rho)\geq0$ if and only if $\cosh\rho\geq\dfrac{2dH+\sqrt{1-4H^2+d^2}}{1-4H^2}$. Let $\rho_1\geq0$ such that $\cosh\rho_1=\dfrac{2dH+\sqrt{1-4H^2+d^2}}{1-4H^2}$, then $f(\rho_1)=0$ and $\rho_1=0$ if and only if $d=-2H$.
\begin{enumerate}
  \item If $d>-2H$, then $\dfrac{-d}{2H}<\cosh\rho_1$. Consequently the function $u$ is nondecreasing for $\rho\geq\rho_1>0$ and has a nonfinite derivative at $\rho_1$.
  \item If $d=-2H$, then $u^\prime(\rho)=\dfrac{2H\sqrt{\cosh\rho-1}\sqrt{1+4\tau^2\tanh^2(\rho/2)}}{\sqrt{(1-4H^2)\cosh\rho+4H^2+1}}$. Therefore the function $u$ is defined for $\rho\geq0$, it has a zero derivative at $0$ and is nondecreasing for $\rho>0$.
  \item If $d<-2H$, then there exist $\rho_0>\rho_1>0$ such that $\dfrac{-d}{2H}=\cosh\rho_0$. Consequently the function $u$ is defined for $\rho\geq\rho_1>0$ with a nonfinite derivative at $\rho_1$, it is nonincreasing for $\rho_1<\rho<\rho_0$, has a zero derivative at $\rho_0$ and it is nondecreasing for $\rho>\rho_0$.
  \item For any $d$ we have $\displaystyle\lim_{\rho\rightarrow+\infty}u(\rho)=+\infty$.
\end{enumerate}
\end{lema}
Next Lemma, is analogous to Lemma \ref{rl1} in the case $H=1/2$. Observe that $f(\rho)=-2d\cosh^2\rho-(1+d^2)$, thus the set $\{\rho,f(\rho)>0\}$ is nonempty if and only if $d<0$

\begin{lema}\label{rl2} Assume $H=1/2$ and $d<0$. Then $f(\rho)\geq0$ if and only if $\cosh^2\rho\geq\frac{1+d^2}{-2d}$. Let $\rho_1\geq0$ such that $\cosh\rho_1=\frac{1+d^2}{-2d}$, then $f(\rho_1)=0$ and $\rho_1=0$ if and only if $d=-1$.
\begin{enumerate}
  \item If $d\in(-1,0)$, then $\frac{-d}{2H}<\cosh\rho_1$. Consequently the function $u$ is nondecreasing for $\rho\geq\rho_1>0$ and has a nonfinite derivative at $\rho_1$.
  \item If $d=-1$, then $u^\prime(\rho)=\frac{1}{\sqrt{2}}\sqrt{(\cosh\rho-1)(1+4\tau^2\tanh^2(\rho/2))}$. Therefore the function $u$ is defined for $\rho\geq0$, it has a zero derivative at $0$ and is nondecreasing for $\rho>0$.
  \item If $d<-1$ there exist $\rho_0>\rho_1>0$ such that $\frac{-d}{2H}=\cosh\rho_0$. Consequently the function $u$ is defined for $\rho\geq\rho_1>0$ with a nonfinite derivative at $\rho_1$, it is nonincreasing for $\rho_1<\rho<\rho_0$, has a zero derivative at $\rho_0$ and it is nondecreasing for $\rho>\rho_0$.
  \item For any d we have $\displaystyle\lim_{\rho\rightarrow+\infty}u(\rho)=+\infty$.
\end{enumerate}
\end{lema}
The proof of Lemma \ref{rl1} and \ref{rl2} is a straightforward computation. As a consequence of Lemma (\ref{rl1}) and Lemma (\ref{rl2}), we have the next results.
\begin{teorema}\label{TR1} (Rotational $H$-surface with $0<H\leq1/2$). Assume $0<H\leq1/2$. there exist a one-parameter family $\mathfrak{H}_d$, $d\in\R$ for $H<1/2$ and $d<0$ for $H=1/2$, of complete rotational $H$-surfaces.
\begin{enumerate}
  \item For $d>-2H$, the surface $\mathfrak{H}_d$ is a properly embedded annulus, symmetric with respect to the slice $\{t=0\}$, the distance between the $"neck"$ and the rotational axis $R=\{(0,0)\times\R\}$ is $arccosh(\dfrac{2dH+\sqrt{1-4H^2+d^2}}{1-4H^2})$ for $H<1/2$ and $arcosh(\dfrac{1+d^2}{-2d})$ for $H=1/2$. See Fig. 1.a
  \item For $d=-2H$, the surface $\mathfrak{H}_{-2H}$ is an entire vertical graph, denoted by $S^H$. Moreover $S^H$ is contained in the halfspace $\{t\geq0\}$ and it is tangent to slice $\mathbb{D}^2\times\{0\}$ at the point $(0,0,0)$. See Fig. 1.b
  \item For $d<-2H$, the surface $\mathfrak{H}_d$ is a properly immersed (and nonembedded) annulus, it is symmetric with respect to slice $\{t=0\}$, the distance between the $"neck"$ and the rotational axis $R$ is $arcosh(\dfrac{2dH+\sqrt{1-4H^2+d^2}}{1-4H^2})$ for $H<1/2$ and $arcosh(\frac{1+d^2}{-2d})$ for $H=1/2$. See Fig. 1.c
  \item In each of the previous case the surface is unbounded in the $t$-coordinate. When $d$ tends to $-2H$ with either $d>-2H$ or $d<-2H$, then the surface $\mathfrak{H}_d$ tends toward the union of $S^H$ and its symmetric with respect to the slice $\{t=0\}$. Furthermore, any rotational $H$-surface with $0<H\leq1/2$ is up to an ambient isometry, a part of a surface of the family $\mathfrak{H}_d$.
\end{enumerate}
\end{teorema}

\begin{proof} The result is a straightforward consequence of Lemma \ref{rl1} and Lemma \ref{rl2}. For $d=-2H$, $\mathfrak{H}_{-2H}$ is the rotational surface generated by the graph of the function $u$.

For $d\neq-2H$, let $\gamma$ be the union of the graph of $u$ join with its symmetric with respect to the slice $\{t=0\}$. Then $\mathfrak{H}_d$ is the rotational surface generated by the curve $\gamma$.
\end{proof}

\begin{figure}[h]
\begin{center}
\psfrag{a}{$a$}
\psfrag{b}{$b$}
\psfrag{c}{$c$}
\psfrag{d}{$1$}
\psfrag{F}{Schematic Figures}
\includegraphics[width=5in]{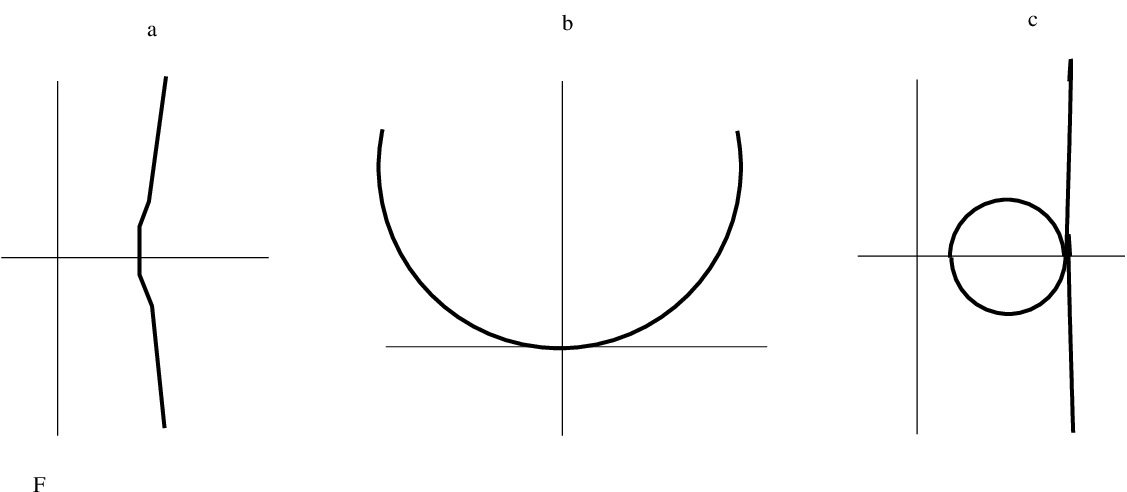}
\end{center}
\end{figure}
\begin{center}
Figure $1.-$ Generating curve for rotational surfaces with $H\leq1/2$.
\end{center}

.\newline

Observe that, $f(\rho)=(1-4^2)\cosh^2\rho-4Hd\cosh\rho-(1+d^2)$, so for $H>1/2$, the set $\{\rho,f(\rho)>0\}$ is nonempty if and only if $d<0$. Furthermore, $f\left(\dfrac{2dH\pm\sqrt{1-4H^2+d^2}}{1-4H^2}\right)=0$, this least equality is possible since $1-4H^2+d^2>0$, this is $d<-\sqrt{4H^2-1}$.

\begin{lema} \label{rl3} Let $H$ and $d$ satisfying $H>1/2$ and $d<-\sqrt{4H^2-1}$. Then, there exist two numbers $0\leq\rho_1<\rho_2$ such that $\cosh\rho_1=\dfrac{2dH+\sqrt{1-4H^2+d^2}}{1-4H^2}$ and $\cosh\rho_2=\dfrac{2dH-\sqrt{1-4H^2+d^2}}{1-4H^2}$. Therefore, $f(\rho)>0$ if and only if $\rho_1<\rho<\rho_2$ and $f(\rho_1)=f(\rho_2)=0$.
\begin{enumerate}
  \item If $d<-2H$, then $\rho_1>0$ and there exist a unique number $\rho_0\in(\rho_1,\rho_2)$ satisfying $g(\rho_0)=0$. Furthermore $g\leq0$ on $[\rho_1,\rho_0)$ and $g\geq0$ on $(\rho_0,\rho_2]$. Consequently, the function $u$ is defined on $[\rho_1,\rho_2]$, has a nonfinite derivative at $\rho_1$ and $\rho_2$, has a zero derivative at $\rho_0$, is nonincreasing on $(\rho_1,\rho_0)$ and nondecreasing on $(\rho_0,\rho_2)$.
  \item If $d=-2H$, then $\rho_1=0$ and $u^\prime(\rho)=\dfrac{2H\sqrt{\cosh\rho-1}\sqrt{1+4\tau^2\tanh^2(\rho/2)}}{\sqrt{(1-4H^2)\cosh\rho+4H^2+1}}$. Consequently, the function $u$ is defined on $[0,\rho_2]$, is nondecreasing, has a zero derivative at $0$ and a nonfinite derivative at $\rho_2$.
  \item If $-2H<d<-\sqrt{4H^2-1}$, then $\rho_1>0$ and $g\leq0$ on $[\rho_1,\rho_2]$. Therefore the function $u$ is defined on $[\rho_1,\rho_2]$, is nondecreasing and has nonfinite derivative at $\rho_1$ and $\rho_2$.
\end{enumerate}
\end{lema}

An immediate consequence of Lemma \ref{rl3} we obtain the next Theorem.

\begin{teorema}\label{TR2} (Rotational surfaces with $H>1/2$) Assume $H>1/2$. There exist a one-parameter family $\mathfrak{D}_d$ of complete rotational $H$-surfaces, $d\leq-\sqrt{4H^2-1}$.
\begin{enumerate}
  \item For $d<-2H$, the surface $\mathfrak{D}_d$ is an immersed (and nonembedded) annulus, invariant by a vertical translation and is contained in the closed region bounded by the two vertical cylinders $\rho=\rho_1$ and $\rho=\rho_2$. Furthermore $\rho_1\rightarrow+\infty$ and $\rho_2\rightarrow+\infty$ when $d\rightarrow-\infty$ and $\rho_1\rightarrow0$ and $\rho_2\rightarrow arcosh\left(\dfrac{4H^2+1}{4H^2-1}\right)$ when $d\rightarrow-2H$. Such surfaces are analogous to the nodoids of Delaunay in $\R^3$. See Fig. 2.a
  \item For $d=-2H$, the surface $\mathfrak{D}_{-2H}$ is an embedded sphere and the maximal distance from the rotational axis is $\rho_2=arcosh\left(\dfrac{4H^2+1}{4H^2-1}\right)$. See Fig. 2.b
  \item For $-2H<d<-\sqrt{4H^2-1}$; the surface $\mathfrak{D}_d$ is an embedded annulus, invariant by a vertical translation and is contained in the closed region bounded by the two vertical cylinders $\rho=\rho_1$ and $\rho=\rho_2$. Furthermore $\rho_1\rightarrow 0$ and $\rho_2\rightarrow arcosh\left(\dfrac{4H^2+1}{4H^2-1}\right)$ when $d\rightarrow-2H$ and both $\rho_1,\rho_2\rightarrow arcosh\left(\dfrac{2H}{\sqrt{4H^2-1}}\right)$ when $d\rightarrow-\sqrt{4H^2-1}$. Moreover $\rho_2\rightarrow arcosh\left(\dfrac{2H}{\sqrt{4H^2-1}}\right)<\rho_2$. Such surfaces are analogous to the undoloids of Delaunay in $\R^3$. See Fig. 2.c
  \item For $d=-\sqrt{4H^2-1}$, the surface $\mathfrak{D}_{-\sqrt{4H^2-1}}$ is the vertical cylinder over the circle with hyperbolic radius $arcosh\left(\dfrac{2H}{\sqrt{4H^2-1}}\right)$.
\end{enumerate}
\end{teorema}

\begin{figure}[h]
\begin{center}
\psfrag{a}{$a$}
\psfrag{b}{$b$}
\psfrag{c}{$c$}
\psfrag{F}{Schematic Figures}
\includegraphics[width=5in]{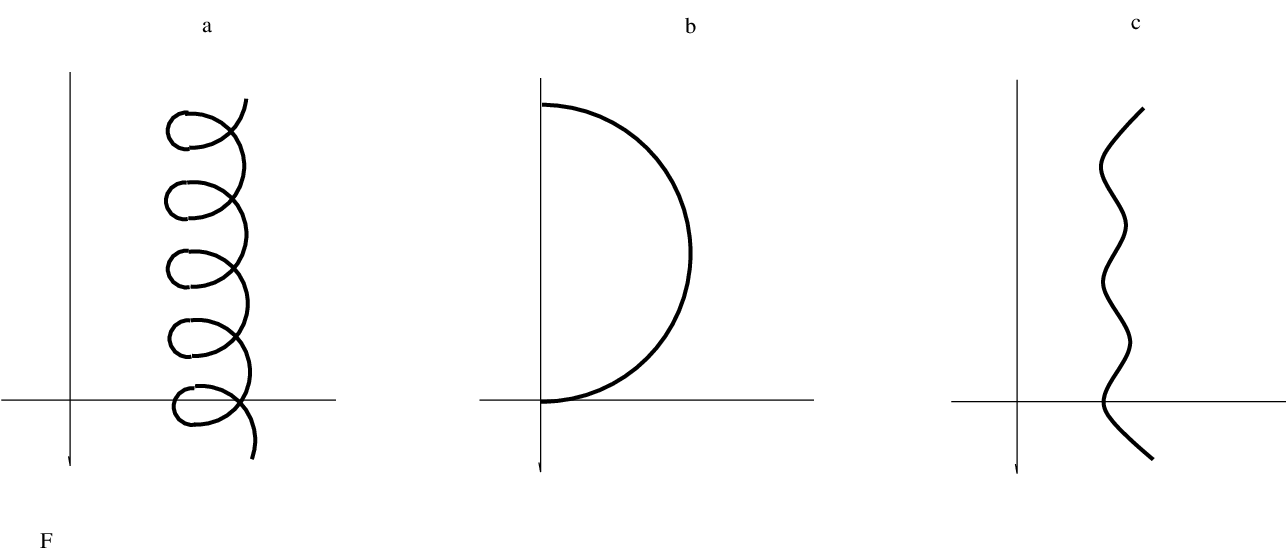}
\end{center}
\end{figure}
\begin{center}
Figure $2.-$ Generating curve for rotational surfaces with $H>1/2$
\end{center}

\section{Parabolic Screw Motions Surfaces in $\widetilde{PSL}_2(\R,\tau)$}

A parabolic screw motion surface in $\widetilde{PSL}_2(\R,\tau)$ is a surface which is invariant by one-parameter of isometries, this one-parameter group of isometries is the composition of parabolic isometries together with vertical translation.

To study this kind of surface we take $M^2=\mathbb{H}^2$ the half space for the model of the Hyperbolic space.\\
The idea is simple, we take a curve in the $xt$ plane and we will apply one-parameter group of parabolic isometries together with translation, to study this kind of surface we take the half plane model, this is $M^2\equiv\mathbb{H}^2$. In this model the parabolic isometries are simply the translation in the $x$ direction.

Let $\alpha(y)=(0.y,u(y))$ a curve in the $xt$ plane and denote by $\Gamma$ the one-parameter group such that the surface $S=\Gamma(\alpha)$ be a parabolic screw motion in $\widetilde{PSL}_2(\R,\tau)$, then $S$ is parameterized by,

  \begin{equation*}
    \varphi(x,y)=(x,y,u(y)+lx)
\end{equation*}
where $l\in\R$

\begin{lema}\label{p1} With the notations above, and denoting by $H$ the mean curvature of $S$, then the function $u$ satisfies
\begin{equation*}
    u(y)=\displaystyle\int\dfrac{(dy-2H)\sqrt{1+(ly-2\tau)^2}}{y\sqrt{1-(dy-2H)^2}}
\end{equation*}
where $d$ is a real number, and $l$ is the pitch.
\end{lema}
\begin{proof} The proof is analogous to the proof of Lemma \ref{r1}.
\end{proof}

\section{Parabolic surfaces in $\widetilde{PSL}_2(\R,\tau)$}

We focus our attention in parabolic surfaces, that is surfaces which are invariant by parabolic isometries, by considering the notation of the screw motion surfaces and making $l\equiv0$ we have the next corollary.

\begin{cor}\label{CP} Denoting by $H$ the mean curvature of $S$, then the function $u$ satisfies
\begin{equation*}
    u(y)=\displaystyle\int\dfrac{(dy-2H)\sqrt{1+4\tau^2}}{y\sqrt{1-(dy-2H)^2}}
\end{equation*}
where $d$ is a real number.
\end{cor}

After a straightforward computation we obtain the next lemma.
\begin{lema}\label{EFP} The solution of the integral is given by
\begin{itemize}
  \item If $H\equiv0$, then
  \begin{equation*}
    u(y)=\sqrt{1+4\tau^2}\arcsin(dy)
  \end{equation*}
  \item If $H=\displaystyle\frac{1}{2}$, then
  \begin{equation*}
   u(y)=\sqrt{1+4\tau^2}\arcsin(dy-1)+\displaystyle\frac{2\sqrt{1+4\tau^2}}{\tan(\frac{\arcsin(cy-1)}{2})+1}
   \end{equation*}
  \item If $H>\displaystyle\frac{1}{2}$, then
   \begin{equation*}
   u(y)=\sqrt{1+4\tau^2}\arcsin(dy-2H)- \displaystyle\frac{4\sqrt{1+4\tau^2}H}{\sqrt{4H^{2}-1}}\arctan\left(\frac{2H\tan(\frac{\arcsin(dy-2H)}{2})+1}{\sqrt{4H^{2}-1}}\right)
   \end{equation*}
\end{itemize}
where $d\in \R$.
\end{lema}

This Lemma gives an immediate examples:

\begin{example} Considering $H\equiv0$, $\tau=-1/2$ and $d=1$, we obtain a parabolic minimal surfaces which is a vertical graph, by considering the rotation by $\pi$ around the $y$ axis we obtain a complete embedded minimal surfaces invariant by parabolic isometries in $\widetilde{PSL}_2(\R,\tau)$.

\begin{center}
  \begin{tabular}{cc}
    \includegraphics[width=2.23in]{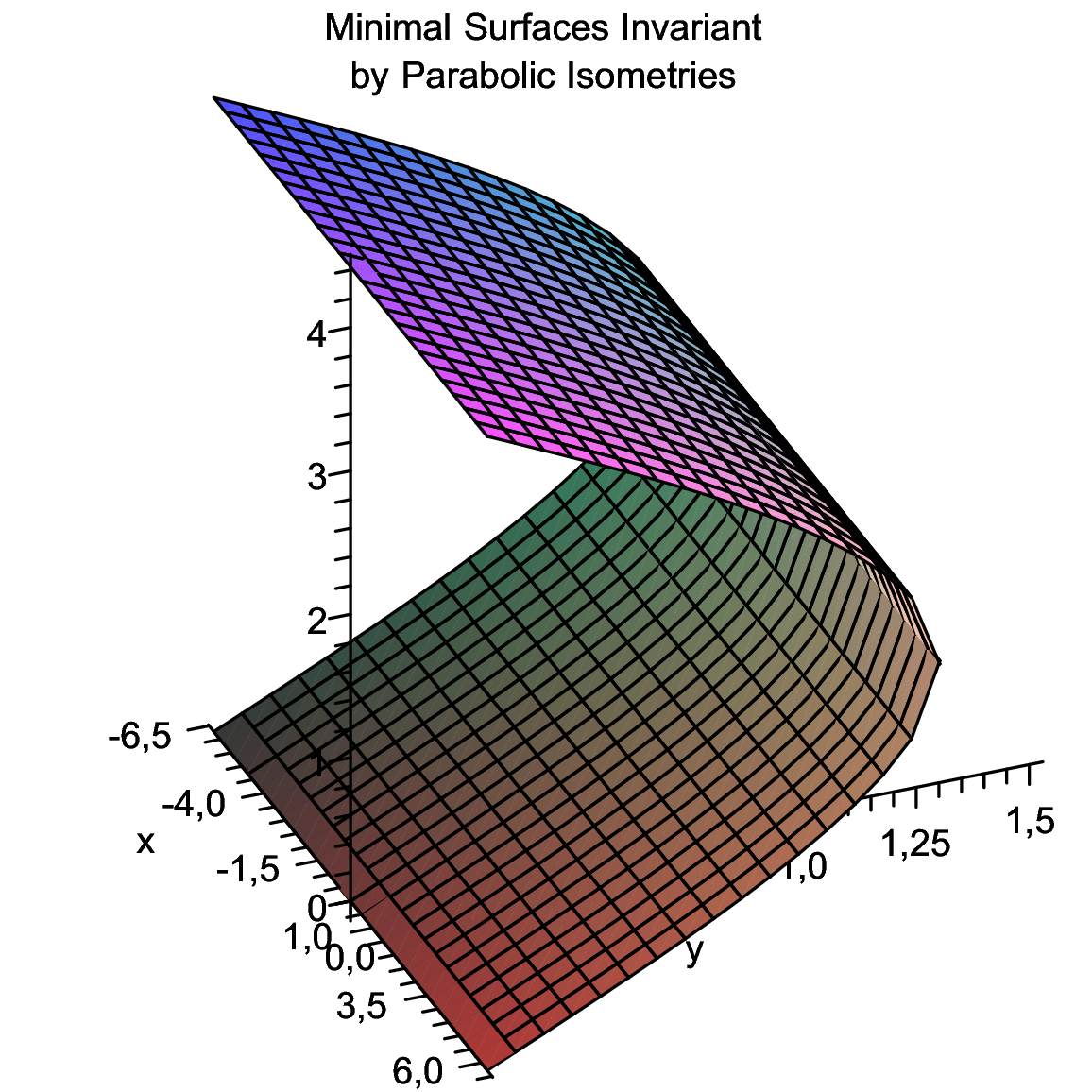} &
    \includegraphics[width=2.23in]{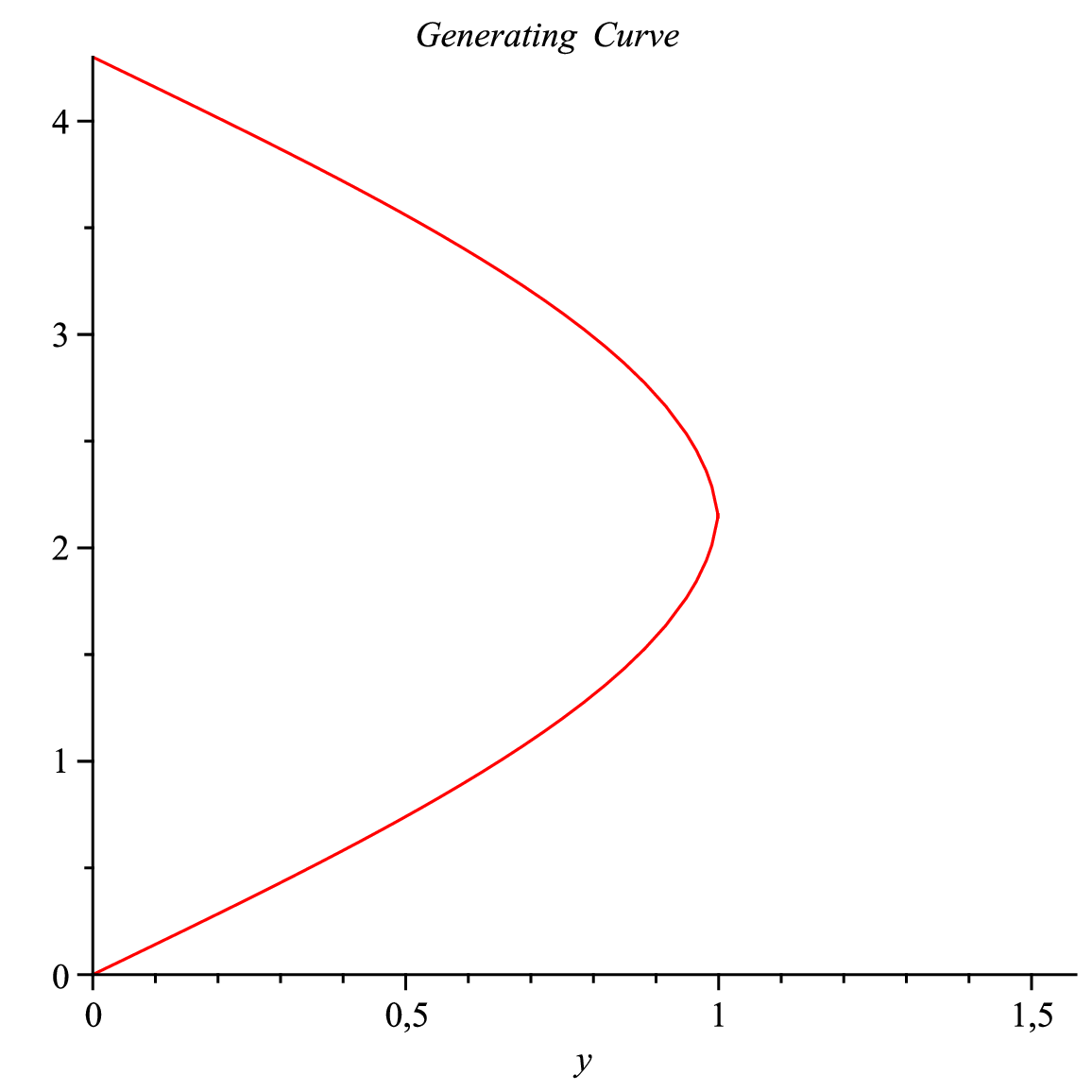}
  \end{tabular}
\end{center}

\end{example}

\begin{example} Considering  $H=1/2$ and $d=1/2$, we obtain
\begin{equation*}
    u(y)=\sqrt{2}\arcsin(dy-1)+\displaystyle\frac{2\sqrt{2}}{\tan(\frac{\arcsin(dy-1)}{2})+1}
\end{equation*}
with Maple's help:

\begin{center}
  \begin{tabular}{cc}
    \includegraphics[width=3in]{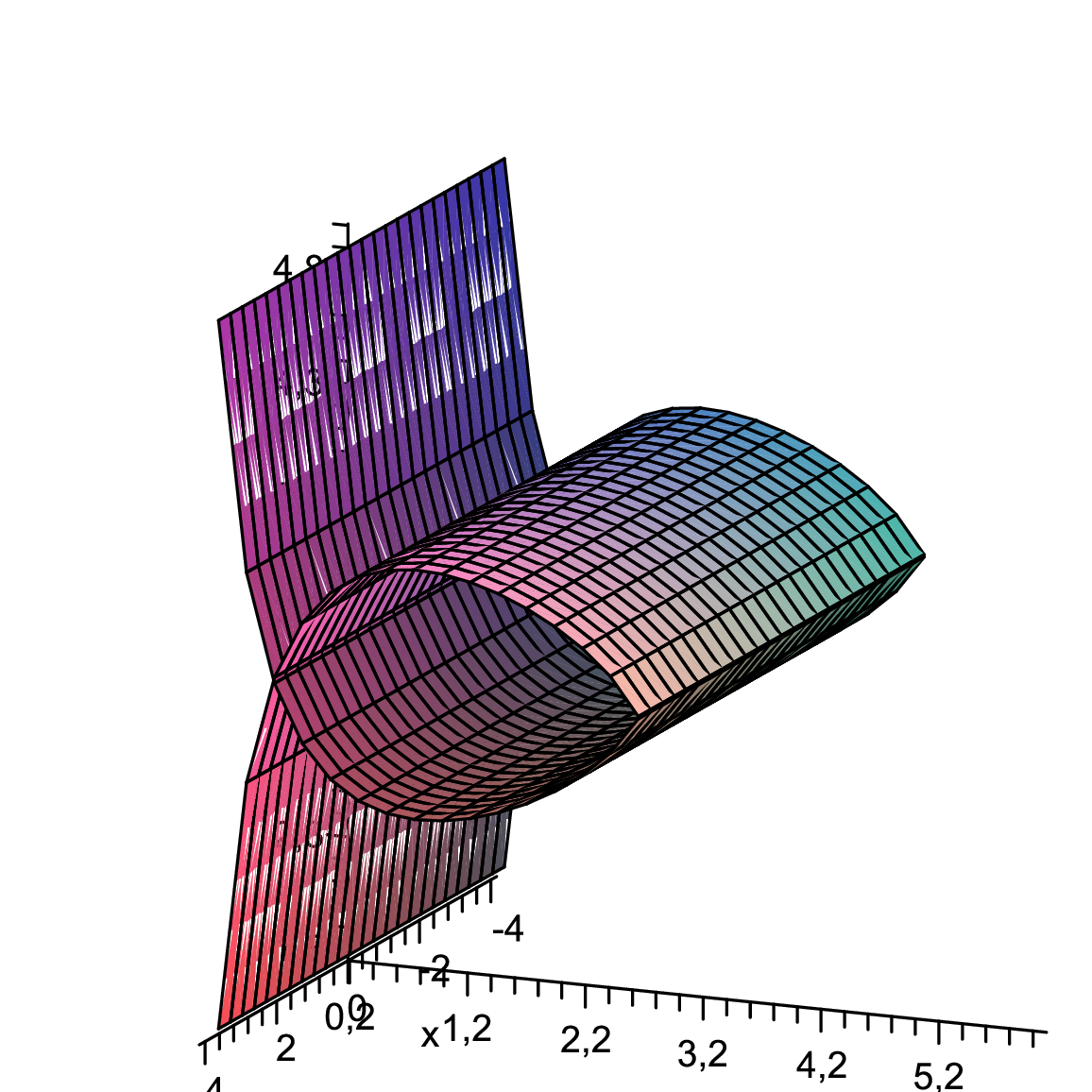} &
    \includegraphics[width=2.23in]{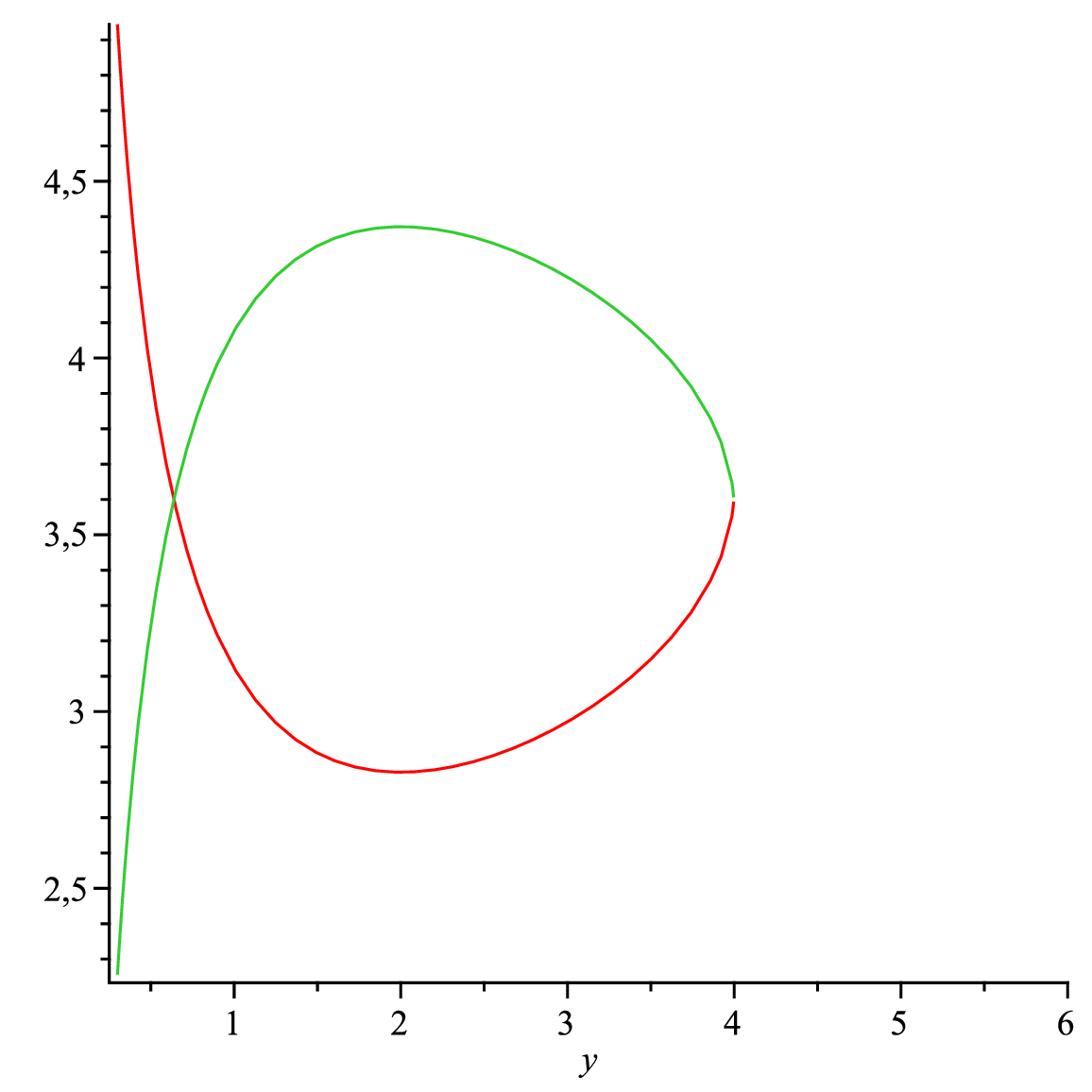}
  \end{tabular}
\end{center}

\end{example}

\begin{example} Finally, we plot a $H=2$ surfaces invariant by Parabolic isometries. Putting $d=8$, $\tau=-1/2$ and $H=2$, we obtain:
\begin{equation*}
    u(y)=\sqrt{2}\arcsin(8y-2H)- \displaystyle\frac{4\sqrt{2}H}{\sqrt{4H^{2}-1}}\arctan \left(\frac{2H\tan(\frac{\arcsin(8y-2H)}{2})+1}{\sqrt{4H^{2}-1}}\right)
\end{equation*}

\begin{center}
  \begin{tabular}{cc}
    \includegraphics[width=4in]{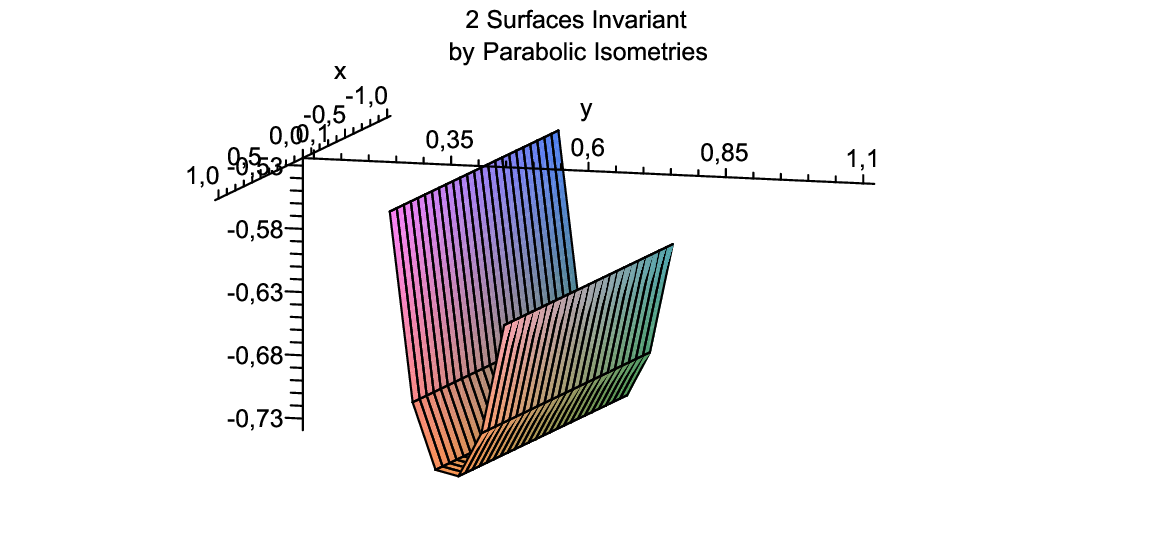} &
    \includegraphics[width=1.8in]{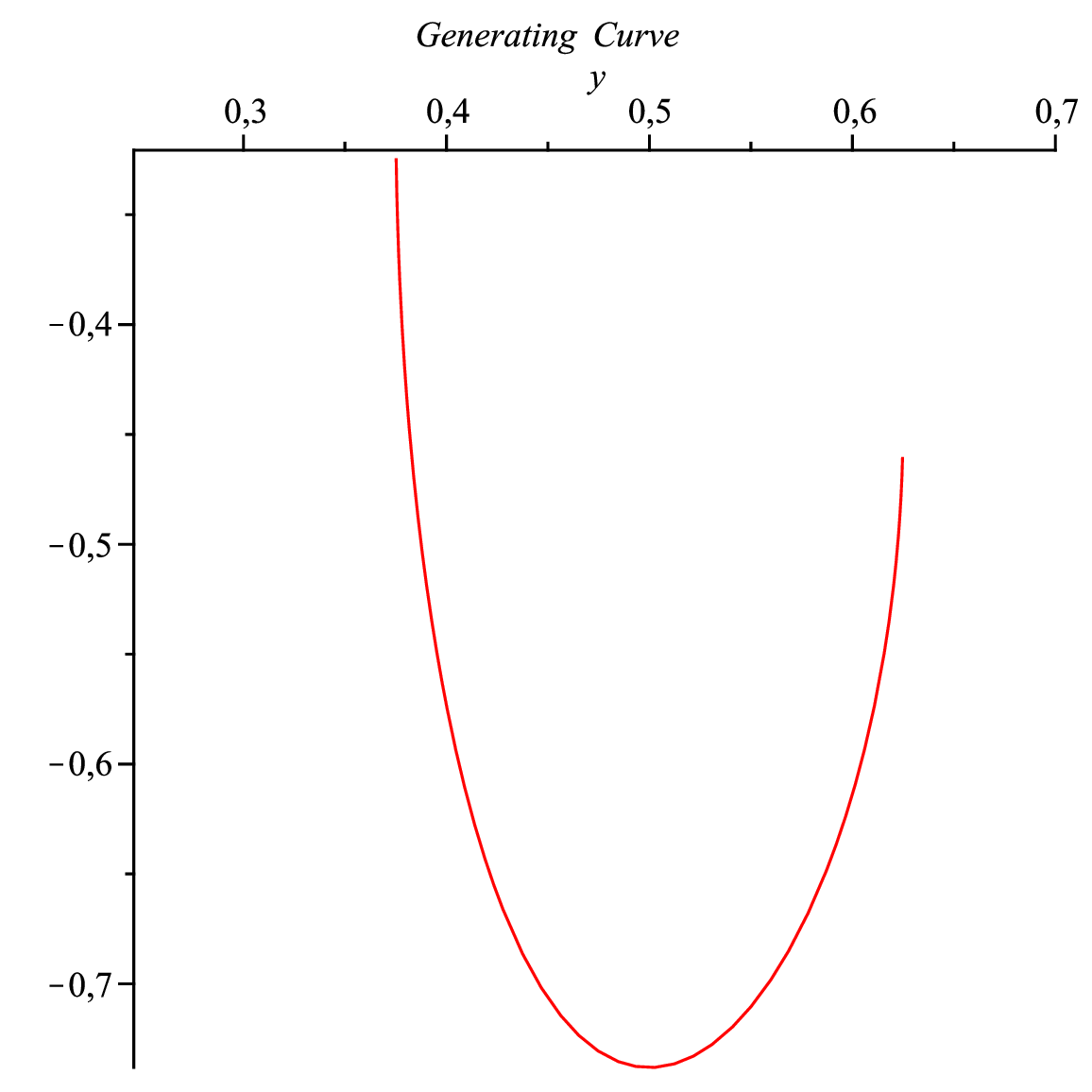}
  \end{tabular}
\end{center}

\end{example}

\section{Surfaces Invariant by Parabolic Isometries in $\widetilde{PSL}_2(\R,\tau)$ with Constant Mean Curvature $H\neq0$}

In this section we describe the behavior of surfaces invariant by parabolic isometries, which have constant mean curvature $H\neq0$. For later use we define the function $g(y)=dy-2H$. Taking into account Formula (\ref{p1}), we obtain the next Lemma

\begin{lema}\label{pl1} Let $H$ be the mean curvature of the surface generated by the Formula (\ref{p1}). Then
\begin{enumerate}
  \item If $d>0$, we have
\begin{itemize}
  \item If $1/2<H$, then $y_1<y<y_2$ where $y_1=\dfrac{2H-1}{d}$ and $y_2=\dfrac{2H+1}{d}$ and there exist a unique number $y_0=\dfrac{2H}{d}\in(y_1,y_2)$ satisfying $g(y_0)$=0. Furthermore $g\leq0$ on $[y_1,y_0)$ and $g\geq0$ on $(y_0,y_1]$. Consequently, the function $h(y)$ is defined on $[y_1,y_2]$, has a nonfinite derivative at $y_1$ and $y_2$, is strictly decreasing on $(y_1,y_0)$ and strictly creasing on $(y_0,y_2)$.
  \item If $0<H<1/2$, then $0<y<y_2$ and there exist a unique number $y_0=\dfrac{2H}{d}\in(0,y_2)$ satisfying $g(y_0)$=0. Furthermore $g\leq0$ on $(0,y_0)$ and $g\geq0$ on $(y_0,y_1]$. Consequently, the function $u(y)$ is defined on $(0,y_2]$, and it is asymptotic to the asymptotic boundary. The function $u$ has a nonfinite derivative at $y_2$, is strictly decreasing on $(0,y_0)$ and strictly creasing on $(y_0,y_2)$.
\end{itemize}
  \item If $d<0$, we have
\begin{itemize}
  \item Here, necessarily $0<H<1/2$. Setting $d=-c$, we have that, $0<y<y_2$, where $y_2=\dfrac{1-2H}{c}$. Consequently, the function $u(y)$ is defined on $(0,y_2]$, has a nonfinite derivative $y_2$, is strictly decreasing on $(0,y_2)$, and $u$ is asymptotic to asymptotic boundary.
\end{itemize}
\end{enumerate}
\end{lema}
As a consequence of Lemma (\ref{pl1}), we have the next results.
\begin{teorema} Let $S$ be the $H$ surface invariant by parabolic isometries immersed into $\widetilde{PSL}_2(\R,\tau)$. Then, there exist a one-parameter family $\mathcal{P}_d$, $d\in\R$ of complete parabolic $H$-surfaces such that,
\begin{enumerate}
  \item For $d>0$, and $H>1/2$ the surface $\mathcal{P}_d$ is immersed (and nonembedded) annulus, invariant by vertical translation, and is contained in the closed region bounded by the vertical cylinders $y=y_1$ and $y=y_2$. See Fig. 3.a
  \item For $d>0$, and $0<H<1/2$ the surface $\mathcal{P}_d$ is a properly immersed (and nonembedded) annulus, it is symmetric with respect to slice $t=0$, the maximum value of $y$ is $y=y_2$. See Fig. 3.b
  \item For $d=-c<0$ and $0<H<1/2$ the surface $\mathcal{P}_d$ is a properly embedded annulus symmetric eith respect to the slice $t=0$, and the maximum value of $y$ is $y=y_2$. See Fig. 3.c
  \item When $d$ tends to $0$, then the surface $\mathcal{P}_d$ tends toward the surface
  \begin{equation*}
    F(y)=\dfrac{-2\sqrt{1+4\tau^2}H\ln(y)}{\sqrt{1-4H^2}}
  \end{equation*}
\end{enumerate}
\end{teorema}

\begin{figure}[h]
\begin{center}
\psfrag{a}{$a$}
\psfrag{b}{$b$}
\psfrag{c}{$c$}
\psfrag{F}{Schematic Figures}
\includegraphics[width=3in]{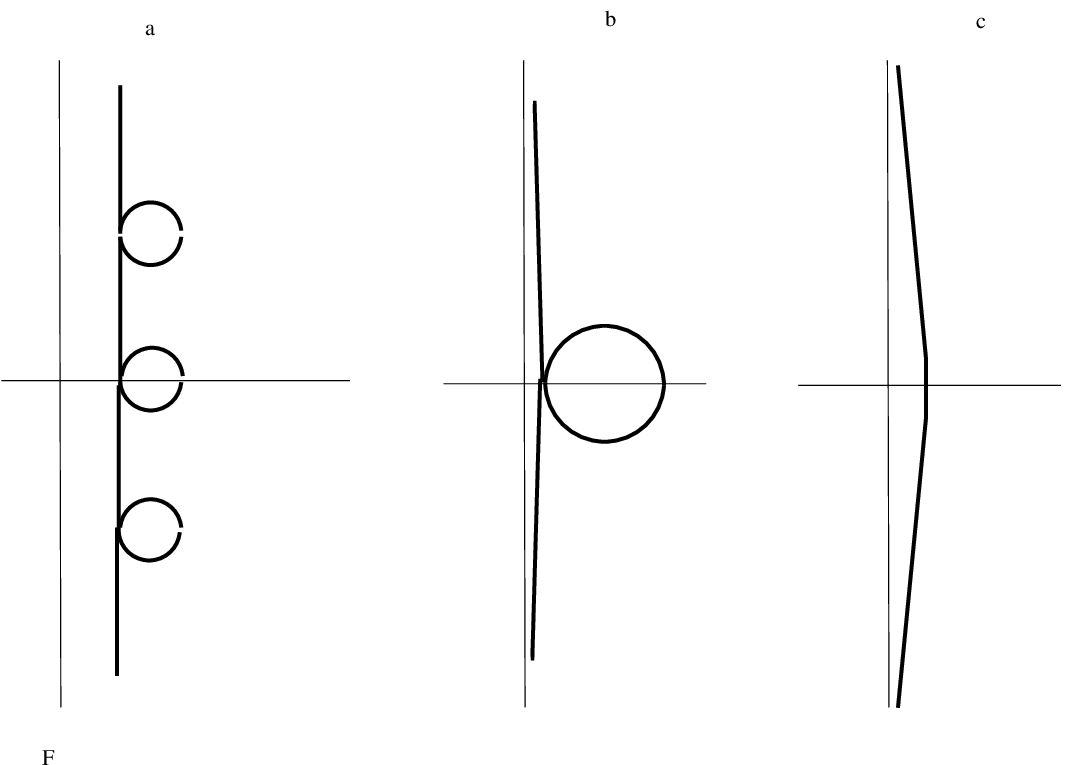}
\end{center}
\end{figure}
\begin{center}
Figure $3.-$ Generating curve for parabolic surfaces with $H\neq1/2$
\end{center}

.\newline

Now, we consider the case $H\equiv1/2$, this is the function $h$ in the formula \ref{CP} become

\begin{equation}\label{Pe1.1}
    h(y)=\sqrt{1+4\tau^2}\int\frac{dy-1}{y}\frac{1}{\sqrt{1-(dy-1)^{2}}}dy
\end{equation}
We denote by $f(y)=1-(dy-1)^2$ and $g(y)=dy-1$, so we obtain the next lemma.

\begin{lema}\label{pl2} By considering the parabolic surface $S$ with constant mean curvature $H=1/2$, we obtain that $d>0$ and the function $h(y)$ is defined for $0<y<y_1=\dfrac{2}{d}$. Furthermore, there exist a number $y_0=\dfrac{1}{d}$ with $0<y_0<y_1\dfrac{2}{d}$ such that $g(y)$ is positive for $0<y<y_0$, $g(y_0)=0$ and $g(y)$ is negative for $y_0<y<y_1$. Consequently the function $h(y)$ is strictly decreasing for $0<y<y_0$, has a horizontal tangent at $y=y_0$ and is strictly increasing for $y_0<y<y_1$. The function $h$ is asymptotic to the asymptotic boundary of $\widetilde{PSL}_{2}(\R,\tau)$
\end{lema}

As a consequence of Lemma \ref{pl2} we have the next result.
\begin{teorema}\label{TP2}  Let $S$ be the $H=1/2$ surface invariant by parabolic isometries immersed into $\widetilde{PSL}_2(\R,\tau)$. Then, there exist a one-parameter family $\mathcal{J}_d$, $d\in\R_+$ of complete parabolic $H$-surfaces such that the surface $\mathcal{J}_d$ is a properly immersed (and nonembedded) annulus, it is symmetric with respect to slice $t=0$, the maximum value of $y$ is $y=y_2$. See Fig. 4
\end{teorema}

\begin{figure}[h]
\begin{center}
\psfrag{a}{$a$}
\psfrag{F}{Schematic Figures}
\includegraphics[width=0.8in]{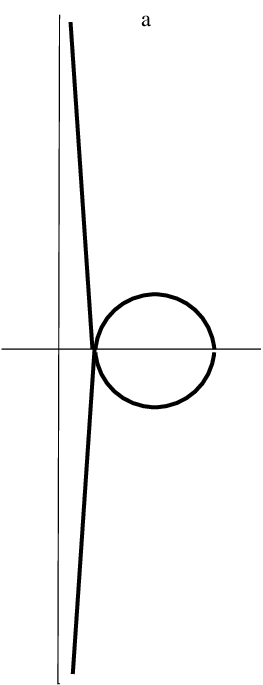}
\end{center}
\end{figure}

\begin{center}
Figure 4..- Generating curve for parabolic surfaces with $H\equiv1/2$
\end{center}

\section{Applications}

In this section we use the study of rotational and parabolic surfaces as well as the examples constructed to prove some general result on graph and multi-graph with constant mean curvature, see Theorem \ref{MT1}.

\begin{prop}\label{AP1} There is no entire graph with constant mean curvature $H$ in $\widetilde{PSL}_2(\R,\tau)$ such that $H>1/2$
\end{prop}
\begin{proof} Suppose that such entire graph exists, without less generality we can suppose that the mean curvature vector field point up. consider the rotational sphere $S$ (here $H>1/2$), given by the Lemma \ref{EFR}, after a vertical translation we can suppose that the sphere $S$ is above of the entire graph and the intersection between this two surface is empty, by considering vertical translation we have a first contact point at the interior of $S$ and the entire graph, by the maximum principle the entire graph is compact, this contradiction complete the proof
\end{proof}

\subsection{Graph and Multi-graph}

Actually the Proposition \ref{AP1} is part of a general theorem see Theorem \ref{d1}.

The classification of simply connected homogeneous manifolds of dimension 3 is well know. Such a manifold has an
isometry group of dimension 3, 4 or 6. When the dimension of the isometry group is 6, we have a space form. When the dimension of the isometry group is 3, the manifold has the geometry of the Lie group $Sol3$.

We will consider the complete homogeneous manifolds $E^{3}(\kappa,\tau)$ whose \newline
isometry groups have dimension 4: such a manifold is a Riemannian fibration over a $2-dimensional$ space form $M^{2}(\kappa)$ that is, $\pi:E^{3}(\kappa,\tau)\longrightarrow M^{2}(\kappa)$ is a Riemannian submersion, where $M^{2}(\kappa)$ is the space form of dimension 2 which has Gauss curvature $\kappa$. If $E^{3}(\kappa,\tau)$ is not compact, then it is topologically  $M^{2}(\kappa)\times\R$, each fiber is diffeomorphic to $\R$ (the real line) and has curvature $\tau$; if $E^{3}(\kappa,\tau)$ is compact, with $\kappa>0$ and $\tau\neq0$, $E^{3}(\kappa,\tau)$ are the compact Berger spheres, each fiber is diffeomorphic to $S^{1}$ (the unit circle).  the tangent unit vector field to the fiber is an unit Killing field which we will denote by $E_{3}$, this vector field will be called the vertical vector field. These manifolds are classified, up to isometry, by the curvature $k$ of the base surface of the fibration and the bundle curvature $\tau$, where $\kappa$ and $\tau$ can be any real numbers satisfying $\kappa\neq 4\tau^{2}$. Namely, these manifolds are
\begin{itemize}
  \item $E^{3}(\kappa,\tau)=\mathbb{D}^{2}\left(\dfrac{1}{\sqrt{-\kappa}}\right)\times\R$, if $\kappa<0$ and $\tau=0$
  \item $E^{3}(\kappa,\tau)=\mathbb{S}^{2}(\sqrt{\kappa})\times\R$, if $\kappa>0$ and $\tau=0$
  \item $E^{3}(\kappa,\tau)=Nil_{3}$ (Heisenberg space) if $\kappa=0$ and $\tau\neq 0$
  \item $E^{3}(\kappa,\tau)=\widetilde{PSL}_{2}(\R)$, if $\kappa<0$ and $\tau\neq 0$
  \item $E^{3}(\kappa,\tau)=\mathbb{S}^{3}_{\tau}$ (Spheres of Berger), if $\kappa>0$ and $\tau\neq 0$
\end{itemize}

Let $M^{n}$ be a Riemannian manifold of dimension $n$, and $\Omega\subset M^{n}$ an open domain in $M^{n}$, such that $\overline{\Omega}$ is compact and $\partial\Omega$ is of class $C^{\infty}$. The Cheeger constant which is denoted by $C(M^{n})$ is given by
\begin{equation*}
    C(M^{n})=\inf_{\Omega}\{\frac{A(\partial\Omega)}{V(\Omega)};\Omega\subset M^n,\overline{\Omega}=compact \}
\end{equation*}
where $A$ is the area function and $V$ the volume function on $M^{n}$.
\begin{obs} We only consider $\kappa\in\{-1,0,1\}$, thus $C(\mathbb{D}^2)\equiv1$.
\end{obs}
\begin{teorema}\label{d1} There is no entire $H$-$graph$ into $\pi:E^3(\kappa,\tau)\longrightarrow M^{2}(\kappa)$ such that $2H>C(M^{2}(\kappa))$, where $C(M^{2}(\kappa))$ is the Cheeger constant of $M^{2}(\kappa)$.
\end{teorema}

This theorem have an important consequence, see Theorem \ref{MT1}

We consider an immersed surface $\Sigma$ with constant mean curvature $H$ in $E^{3}(\kappa,\tau)$. Let $N$ be the normal unit vector field along $\Sigma$, Let $\nu=g(N,E_{3})$ be the angle function on $\Sigma$.
\begin{defi} $\Sigma$ is said $Multi-graph$ if the function $\nu$ satisfy, either $\nu>0$ or $\nu<0$ on $\Sigma$.
\end{defi}

\begin{teorema}\label{MT1} Let $\Sigma$ be a complete $H$ surfaces immersed into $E^3(\kappa,\tau)$, such that $\nu$ does not change of sign, that is either $\nu\geq0$ or $\nu\leq0$. Then if $2H>C(M^{2}(\kappa))$, where $C(M^{2}(\kappa))$ is the Cheeger constant of $M^{2}(\kappa)$, we have that $\nu\equiv0$, ie. $\Sigma$ is a vertical cylinder.

In particular, there is no immersed complete $H$ $multi-graph$ $\Sigma$ in $E(\kappa,\tau)$ such that $2H>C(M^{2}(\kappa))$.
\end{teorema}
 The details of the proof are in \cite{CP}, other related references are \cite{RRT}, \cite{H-E}, \cite{HRS} \cite{Chg}, \cite{CBE}, \cite{CB}, and \cite{D}.

\subsection{Asymptotic Behavior of Rotational Surface}

In the paper, A halfspace theorem for mean curvature $H=\frac{1}{2}$ surfaces in $\mathbb{H}^2\times\R$ (see \cite[Theorem 1]{BS}), the authors studied the asymptotic behavior for $H=1/2$ rotational surfaces immersed into $\mathbb{H}^2\times\R$, and prove a halfspace theorem, since the behavior of $H=1/2$ rotational surfaces immersed into $\widetilde{PSL}_2(\R,\tau)$ is similar rotational surfaces in $\mathbb{H}^2\times\R$ (see Theorem \ref{TR1}), it is natural ask for a halfspace theorem for mean curvature $H=1/2$ surfaces in $\mathbb{H}^2\times\R$.

Denote by $\alpha=-d$, then for any $\alpha\in\R_+$, there exist a rotational surface $\mathfrak{H}_{\alpha}$ of constant mean curvature $H=\frac{1}{2}$. For $\alpha\neq1$, the surface $\mathfrak{H}_{\alpha}$ has two vertical ends ( where a vertical end is a topological annulus, with no asymptotic point at finite height) that are vertical graph over the exterior of a disk $D_{\alpha}$.

Up to vertical translation, one can  assume that $\mathfrak{H}_{\alpha}$ is symmetric with respect to the horizontal plane $t=0$. For $\alpha=1$, the surface $\mathfrak{H}_1$ has only one end, it is a graph over $\mathbb{D}^2$ and it is denoted by $S$.

For any $\alpha>1$ the surface $\mathfrak{H}_{\alpha}$ is not embedded. The self-intersection set is a horizontal circle on the plane $t=0$. Denote by $\rho_{\alpha}$ the radius of the intersection circle. For $\alpha<1$ the surface $\mathfrak{H}_{\alpha}$ is embedded.

For any $\alpha\in\R_+$, let $u_{\alpha}:\mathbb{D}^2\times\{0\}\backslash D_{\alpha}\longrightarrow\R$ be the function such that the end of the surface $\mathfrak{H}_{\alpha}$ is the vertical graph of $u_{\alpha}$. The asymptotic behavior has the following form: $u_{\alpha}(\rho)\simeq \dfrac{\sqrt{1+4\tau^2}}{\sqrt{\alpha}}\exp^{\left(\frac{\rho}{2}\right)}$, $\rho\rightarrow\infty$, where $\rho$ is the hyperbolic distance from the origin. The positive number $\dfrac{\sqrt{1+4\tau^2}}{\sqrt{\alpha}}\in\R_+$ is called the growth of the end. (see \cite{CP}).

Following the ideas of Barbara Nelli e Ricardo Sa Earp, we have the next theorem.

\begin{teorema}\label{THS1} Let $S$ be a simply connected rotational surface in $\widetilde{PSL}_2(\R,\tau)$, with constant mean curvature $H=\frac{1}{2}$. Let $\Sigma$ be a complete surface with constant mean curvature $H=\frac{1}{2}$, different from a rotational simply connected one. Then, $\Sigma$ cannot be properly immersed in the mean convex side of $S$
\end{teorema}
The proof is analogous to this in $\mathbb{H}^2\times\R$, for more details see \cite{CP}.

\end {document}